\documentstyle[12pt]{article}
\def\a{\alpha}            \def\b{\beta}
\def\g{\gamma}            
\def\d{\delta}            \def\D{\Delta}
\def\th{\theta}           
\def\vth{\vartheta}       
         \def\vep{\varepsilon}
\def\vph{\varphi}
\def\om{\omega}           
               
\def\sg{\sigma}            
              
\def\lm{\lambda}          
         \def\Up{\Upsilon}
             
\def\t{\tau}

\def\tl{\tilde}           
\def\ul{\underline}       \def\ol{\overline}

         \def\cK{{\cal K}}         
                  
                  \def\cR{{\cal R}}

\newcommand{\id}{\mathop{\rm id}\nolimits}

\def\lb{\label}                 \def\ot{\otimes}
\def\cd{\cdots}                 \def\ld{\ldots}

\def\la{\leftarrow}             \def\ra{\rightarrow}

          \def\sms{\small}
            \def\ls{\large}
                 
\def\uD{\underline{\Delta}}
\def\Dp{\Delta_{+}}
\def\uDp{\underline{\Delta}_{+}}

\def\r#1{(\ref{#1})}

\def\NN{{\rm N}\!\!\!\!\!{\rm N}}
\def\ZZ{Z\!\!\!\!\!\!Z}

\def\Z2{Z\!\!\!\!\!\!Z_{\,2}}
\def\CC{{\bf C}\!\!\!{\rm l}}

\newcommand{\bn}{\begin{equation}}
\newcommand{\ed}{\end{equation}}
\newcommand{\bneqn}{\begin{eqnarray}}
\newcommand{\edeqn}{\end{eqnarray}}
\newcommand{\bnth}{\begin{theorem}}
\newcommand{\edth}{\end{theorem}}
\newcommand{\bnpr}{\begin{proposition}}
\newcommand{\edpr}{\end{proposition}}
\newcommand{\bnlm}{\begin{lemma}}
\newcommand{\edlm}{\end{lemma}}
\newcommand{\bndef}{\begin{definition}}
\newcommand{\eddef}{\end{definition}}
\newcommand{\bntb}{\begin{tabular}}
\newcommand{\edtb}{\end{tabular}}

\def\nin{\noindent}
\newtheorem{definition}{Definition}[section]
\newtheorem{lemma}{Lemma}[section]
\newtheorem{proposition}{Proposition}[section]
\newtheorem{theorem}{Theorem}[section]

\def\r#1{(\ref{#1})}
\def\DD{\Delta^{(D)}_{\;q}}
\def\SD{S^{(D)}_{\;q}}
\def\RD{{\cal R}^{(D)}}

\def\kr{\bar{{\cal R}}}
\def\stackreb#1#2{\mathrel{\mathop{#2}\limits_{#1}}}
\normalbaselines
\begin{document}
\begin{center}
\vspace*{1cm}

{\LARGE{\bf Quantum Affine (Super)Algebras\\[3pt]
$U_q(A_{1}^{(1)})$ and $U_q(C(2)^{(2)})$}}\\

\vskip 1.2cm
{\large{\bf S.M. Khoroshkin$^{1}$, J. Lukierski$^{2}$ and
V.N. Tolstoy$^{3}$}}

\vskip 0.5cm
$^1$ Institute of Theoretical \& and Experimental Physics \\
117259 Moscow, Russia (e-mail: khor.@heron.itep.ru)\\

\vskip 0.2cm
$^{2}$Institute of Theoretical Physics, University of Wroc\l{aw}\\
50-204 Wroc\l{aw} \& Poland (e-mail: lukier@ift.uni. wroc.pl)\\

\vskip 0.2cm
$^3$ Institute of Nuclear Physics, Moscow State University \\
119899 Moscow, Russia (e-mail: tolstoy@nucl-th.npi.msu.su)
\end{center}
\date{}

\vskip 0.6cm
\begin{abstract}
We show that the quantum affine algebra $U_{q}(A_{1}^{(1)})$ and
the quantum affine superalgebra $U_{q}(C(2)^{(2)})$ admit a unified
description. The difference between them consists in the phase
factor which is equal to $1$ for $U_{q}(A_{1}^{(1)})$ and  it is equal
to $-1$ for $U_{q}(C(2)^{(2)})$. We present such a description for
the actions of the braid group, for the construction of Cartan-Weyl
generators and their commutation relations, as well for
the extremal projector and the universal $R$-matrix.
We give also a unified description for the 'new
realizations' of these algebras together with explicit calculations
of corresponding $R$-matrices.
\end{abstract}

\setcounter{equation}{0}
\section{Introduction}
Among variety of all affine Lie (super)algebras\footnote{We conclude
the prefix "super" in brackets to stress that the Lie (super)algebras
include both the Lie algebras and the own Lie superalgebras.} (both
quantized and non-quantized) the affine (super)algebras of rank 2
play same key role. In the first place, all affine series of the
type $A(n|m)^{(1)}$, $B(n|m)^{(1)}$, $C(n)^{(1)}$, $D(n|m)^{(1)}$,
$A(2n|2m\!-\!1)^{(2)}$, $A(2n\!-\!1|2m\!-\!1)^{(2)}$,$C(n)^{(2)}$,
$D(n|m)^{(2)}$ and $A(2n|2m)^{(4)}$ are started from the affine
(super)algebras of rank 2. Secondly, the contragredient Lie
(super)algebras of rank 2 are basic structural blocks of any
affine (super)algebras of arbitrary rank. This fact permits,
for example, to reduce the proof of basic theorems for extremal
projector and the universal $R$-matrix to the proof of such theorems
for the (super)algebras of rank 2 (see Refs. \cite{AST},
\cite{T1}--\cite{TK}, \cite{KT1}--\cite{KT5}).
Moreover the representation theory of the affine
(super)algebras (both quantized and non-quantized) contains some
typical elements of the representation theory of the affine
(super)algebras of rank 2. Besides, in  applications of the affine
(super)algebras, first of all the affine (super)algebras of rank 2
are used by virtue of their simplicity.

In this paper we give detailed description of the quantum untwisted
affine algebra $U_{q}(A_{1}^{(1)})$ ($\simeq U_{q}(\hat{sl}(2))$)
and the quantum twisted affine superalgebra $U_{q}(C(2)^{(2)})$
($\simeq U_{q}(\widehat{osp}(2|2))^{(2)}$).
Moreover our goal is to show that these quantum (super)algebras
are described in unified way.
Namely, we present in unified way their defining relations and actions
of the braid group associated with the Weyl group, the construction of
the Cartan-Weyl bases, the complete list of all permutation relations
of the Cartan-Weyl generators corresponding to all root vectors and
finally the unified formula for their extremal projector and universal
$R$-matrix. We extend also the unified description to so called
'new realizations' of the algebras.
Here we present a unified description of the universal $R$-matrices for
corresponding Hopf structures in a multiplicative form as well as
in the form of contour integrals. Difference between both
considered here quantum (super)algebras is only determined by a phase
factor which is equal to $1$ for $U_{q}(A_{1}^{(1)})$ and it is equal
to $-1$ for $U_{q}(C(2)^{(2)})$.
This situation is similar to the finite-dimensional case. Namely, in
the paper \cite{KT1} it was shown that all quantum (super)algebras
$U_q(g)$, where $g$ are the finite-dimensional contragredient Lie
(super)algebras of rank 2, are divided into three classes. Each such
class is characterized by the same Dynkin diagram and the same reduced
root system, provided that we neglect 'colour' of the roots, and all
(super)algebras of the same class have the unified defining relations,
unified construction and  properties of the Cartan-Weyl basis and
unified formula for the universal $R$-matrix. Difference between
the (super)algebras of the same class is determined by some phase
factors which takes values $\pm1$ depending on the colour of the nodes
of their Dynkin diagram.
Concerning the Cartan Weyl bases for the quantum affine algebra
$U_{q}(A_{1}^{(1)})$ and quantum affine superalgebra
$U_{q}(C(2)^{(2)})$ it should be noted that certain results presented
here can be founded in literature separately for the each case
(e.g., see Refs. \cite{Be1}, \cite{D2}, \cite{KT2}--\cite{KT5},
and \cite{YZ}).

Basic information about the (super)algebras $A_{1}^{(1)}$ and
$C(2)^{(2)}$ is presented in the tables 1a and 1b (see Refs. \cite{K1},
\cite{K2}, \cite{VdL}). In the table 1a there are listed
the standard and symmetric Cartan matrices $A$ and $A^{sym}$,
the corresponding extended symmetric matrices $\bar{A}^{sym}$ and
their inverses $(\bar{A}^{sym})^{-1}$, as well as the sets of odd simple
roots (odd roots), the Dynkin diagrams (diagram), and the dimensions of
these (super)algebras (dim).
We remind some elementary definitions of the colour of the roots:
\begin{itemize}
\item
All even roots are called white roots. A white root is
pictured by the white node
\mbox{\begin{picture}(10,10)\put(4,3){\circle{8}}\end{picture}}.
\item
An odd root $\g$ is called a grey root if $2\g$ is not any root.
This root is pictured by the grey node $\ot$.
\item
An odd root $\g$ is called a dark root if $2\g$ is a root.
This root is pictured by the dark node
\mbox{\begin{picture}(10,10)\put(4,3){\circle*{8}}\end{picture}}.
\end{itemize}
We also remind the definition of the reduced system of the positive
root system $\D_+$ for any contragredient (super)algebras of finite
growth.
\begin{itemize}
\item
The system $\ul{\D}_+$ is called the reduced system if it is defined
by the following way:
$\ul{\D}_+\!=\D_+\!\backslash\{2\g\in \D_+|\,\g\;{\rm is\; odd}\}$.
That is the reduced system $\ul{\D}_+$ is obtained from the total
system $\D_+$ by removing of all doubled roots $2\g$ where $\g$ is
a dark odd root.
\end{itemize}
The total and reduced root systems of the (super)algebras $A_{1}^{(1)}$
and $C(2)^{(2)}$ are listed in the table 1b. It is convenient to
present the total root systems  $\D\!=\!\D_+\!\bigcup(-\D_+)$ and reduced
 root systems $\uD\!=\!\uDp\!\bigcup(-\uDp)$
 by the pictures: Figs. 1, 2a, 2b.
Comparing Fig. 1 and Fig. 2b we see that the reduced root systems of
$A_{1}^{(1)}$ and $C(2)^{(2)}$ coincide if we neglect colour of the
roots.

\vskip 15pt
\begin{center}
{\bf Table 1a}

\vskip 10pt
\nin
{\footnotesize
{\renewcommand{\arraystretch}{0}
\bntb{|lccccc|}
\hline
\rule{0pt}{8pt}&&&&&\\
\strut $g(A,\Up)$ & $A=A^{sym}$ & $\bar{A}^{sym}$
&$(\bar{A}^{sym})^{-1}$ & odd & diagram \\[20pt]
\hline
\rule{0pt}{8pt}&&&&&\\
\strut $A_{1}^{(1)}$ & $\left(\hspace{-3mm}
\begin{array}{rr}\strut 2&-2\\
\strut -2&2 \end{array}\right)$
& $\left(\begin{array}{rrr}\strut 0&1&0\\
\strut 1&2&-2\\
\strut 0&-2&2\end{array}\right)$ &
$\left(\begin{array}{rrr}\strut 0&1&1\\
\strut 1&0&0\\
\strut 1&0&\frac{1}{2}\end{array}\right)$ & $\emptyset$
& \makebox{\begin{picture}(10,10)\thicklines
\put(-10,3){\circle{8}}
\put(-20,12){\footnotesize$\d\!-\!\a$}
\put(-7,5){\line(1,0){19}}
\put(-7,1){\line(1,0){19}}
\put(16,3){\circle{8}}
\put(14,12){\footnotesize $\a$}\end{picture}}\\[25pt]
\strut $C(2)^{(2)}$ & $\left(\hspace{-3mm}
\begin{array}{rr}\strut 2&-2\\
\strut -2&2 \end{array}\right)$
& $\left(\begin{array}{rrr}\strut 0&1&0\\
\strut 1&2&-2\\
\strut 0&-2&2\end{array}\right)$ &
$\left(\begin{array}{rrr}\strut 0&1&1\\
\strut 1&0&0\\
\strut 1&0&\frac{1}{2}\end{array}\right)$ & $\{\d\!-\!\a,a\}$
& \makebox{\begin{picture}(10,10)\thicklines
\put(-10,3){\circle*{8}}
\put(-20,12){\footnotesize$\d\!-\!\a$}
\put(-7,5){\line(1,0){19}}
\put(-7,1){\line(1,0){19}}
\put(16,3){\circle*{8}}
\put(14,12){\footnotesize $\a$}\end{picture}}\\[25pt]
\hline
\edtb}}

\vskip 15pt
{\bf Table 1b}

\vskip 10pt
\nin
{\footnotesize
{\renewcommand{\arraystretch}{0}
\bntb{|lcc|}
\hline
\rule{0pt}{8pt}&&\\
\strut $g(A,\Up)$\hspace{1.1cm} & $\Dp$ \hspace{1cm} & $\uDp$
\\[14pt]
\hline
\rule{0pt}{8pt}&&\\
\strut $A_{1}^{(1)}$ \hspace{1.1cm}&
$\{\a,\,n\d\!\pm\!\a,\,n\d\,|\,n\in\NN\}$
\hspace{1cm} & $\{\a,\,n\d\!\pm\!\a,\,n\d\,|\,n\in\NN\}$
\\[20pt]
\strut $C(2)^{(2)}$ \hspace{1.1cm}&
$\{\a,\,2\a,\,n\d\!\pm\!\a,\,2n\d\!\pm\!2\a,\,n\d\,|\,n\in\NN\}$
\hspace{1cm} & $\{\a,\,n\d\!\pm\!\a,\,n\d\,|\,n\in\NN\}$
\\[20pt]
\hline
\edtb}}
\end{center}
\vskip 10pt
\begin{center}
\begin{picture}(450,100)
\put(0,72){\ldots}
\put(24,78){\footnotesize$-\!4\d\!+\!\alpha$}
\put(44,72){\circle{5}}
\put(68,78){\footnotesize$-\!3\d\!+\!\alpha$}
\put(88,72){\circle{5}}
\put(112,78){\footnotesize$-\!2\d\!+\!\alpha$}
\put(132,72){\circle{5}}
\put(161,78){\footnotesize$-\!\d\!+\!\alpha$}
\put(176,72){\circle{5}}
\put(218,78){\footnotesize$\alpha$}
\put(220,72){\circle{5}}
\put(255,78){\footnotesize$\d\!+\!\alpha$}
\put(264,72){\circle{5}}
\put(294,78){\footnotesize$2\d\!+\!\alpha$}
\put(308,72){\circle{5}}
\put(338,78){\footnotesize$3\d\!+\!\alpha$}
\put(352,72){\circle{5}}
\put(382,78){\footnotesize$4\d\!+\!\alpha$}
\put(396,72){\circle{5}}
\put(420,72){\ldots}
\put(0,50){\ldots}
\put(35,56){\footnotesize$-\!4\d$}
\put(44,50){\circle{5}}
\put(79,56){\footnotesize$-\!3\d$}
\put(88,50){\circle{5}}
\put(123,56){\footnotesize$-\!2\d$}
\put(132,50){\circle{5}}
\put(169,56){\footnotesize$-\!\d$}
\put(176,50){\circle{5}}
\put(220,50){\vector(0,1){19}}
\put(220,50){\vector(2,-1){42}}
\put(262,56){\footnotesize$\d$}
\put(264,50){\circle{5}}
\put(303,56){\footnotesize$2\d$}
\put(308,50){\circle{5}}
\put(347,56){\footnotesize$3\d$}
\put(352,50){\circle{5}}
\put(391,56){\footnotesize$4\d$}
\put(396,50){\circle{5}}
\put(4200,50){\ldots}
\put(0,28){\ldots}
\put(24,33){\footnotesize$-\!4\d\!-\!\alpha$}
\put(44,28){\circle{5}}
\put(68,33){\footnotesize$-\!3\d\!-\!\alpha$}
\put(88,28){\circle{5}}
\put(112,33){\footnotesize$-\!2\d\!-\!\alpha$}
\put(132,28){\circle{5}}
\put(161,33){\footnotesize$-\!\d\!-\!\alpha$}
\put(176,28){\circle{5}}
\put(213,33){\footnotesize$-\!\alpha$}
\put(220,28){\circle{5}}
\put(255,33){\footnotesize$\d\!-\!\alpha$}
\put(264,28){\circle{5}}
\put(294,33){\footnotesize$2\d\!-\!\alpha$}
\put(308,28){\circle{5}}
\put(338,33){\footnotesize$3\d\!-\!\alpha$}
\put(352,28){\circle{5}}
\put(382,33){\footnotesize$4\d\!-\!\alpha$}
\put(396,28){\circle{5}}
\put(420,28){\ldots}
\end{picture}
\end{center}

\vskip -20pt
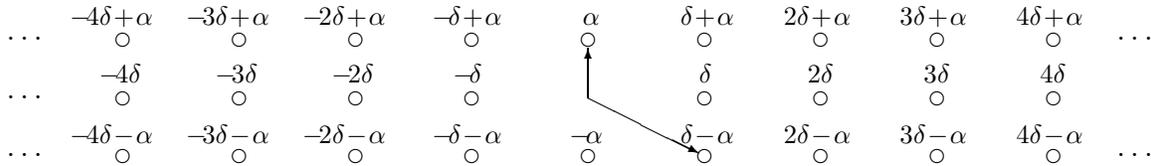
\centerline{\footnotesize Fig. 1. The total and
reduced root system ($\D=\uD$) of $A_{1}^{(1)} (\simeq
\widehat{sl}_{2})$}.

\vskip 10pt
\begin{center}
\begin{picture}(450,120)
\put(0,94){\ldots}
\put(24,101){\footnotesize$-\!4\d\!+\!2\alpha$}
\put(44,94){\circle{5}}
\put(112,101){\footnotesize$-\!2\d\!+\!2\alpha$}
\put(132,94){\circle{5}}
\put(215,101){\footnotesize$2\alpha$}
\put(220,94){\circle{5}}
\put(294,101){\footnotesize$2\d\!+\!2\alpha$}
\put(308,94){\circle{5}}
\put(382,101){\footnotesize$4\d\!+\!2\alpha$}
\put(396,94){\circle{5}}
\put(420,94){\ldots}
\put(0,72){\ldots}
\put(24,78){\footnotesize$-\!4\d\!+\!\alpha$}
\put(44,72){\circle*{5}}
\put(68,78){\footnotesize$-\!3\d\!+\!\alpha$}
\put(88,72){\circle*{5}}
\put(112,78){\footnotesize$-\!2\d\!+\!\alpha$}
\put(132,72){\circle*{5}}
\put(161,78){\footnotesize$-\!\d\!+\!\alpha$}
\put(176,72){\circle*{5}}
\put(218,78){\footnotesize$\alpha$}
\put(220,72){\circle*{5}}
\put(255,78){\footnotesize$\d\!+\!\alpha$}
\put(264,72){\circle*{5}}
\put(294,78){\footnotesize$2\d\!+\!\alpha$}
\put(308,72){\circle*{5}}
\put(338,78){\footnotesize$3\d\!+\!\alpha$}
\put(352,72){\circle*{5}}
\put(382,78){\footnotesize$4\d\!+\!\alpha$}
\put(396,72){\circle*{5}}
\put(420,72){\ldots}
\put(0,50){\ldots}
\put(35,56){\footnotesize$-\!4\d$}
\put(44,50){\circle{5}}
\put(79,56){\footnotesize$-\!3\d$}
\put(88,50){\circle{5}}
\put(123,56){\footnotesize$-\!2\d$}
\put(132,50){\circle{5}}
\put(169,56){\footnotesize$-\!\d$}
\put(176,50){\circle{5}}
\put(220,50){\vector(0,1){19}}
\put(220,50){\vector(2,-1){42}}
\put(262,56){\footnotesize$\d$}
\put(264,50){\circle{5}}
\put(303,56){\footnotesize$2\d$}
\put(308,50){\circle{5}}
\put(347,56){\footnotesize$3\d$}
\put(352,50){\circle{5}}
\put(391,56){\footnotesize$4\d$}
\put(396,50){\circle{5}}
\put(4200,50){\ldots}
\put(0,28){\ldots}
\put(24,33){\footnotesize$-\!4\d\!-\!\alpha$}
\put(44,28){\circle*{5}}
\put(68,33){\footnotesize$-\!3\d\!-\!\alpha$}
\put(88,28){\circle*{5}}
\put(112,33){\footnotesize$-\!2\d\!-\!\alpha$}
\put(132,28){\circle*{5}}
\put(161,33){\footnotesize$-\!\d\!-\!\alpha$}
\put(176,28){\circle*{5}}
\put(213,33){\footnotesize$-\!\alpha$}
\put(220,28){\circle*{5}}
\put(255,33){\footnotesize$\d\!-\!\alpha$}
\put(264,28){\circle*{5}}
\put(294,33){\footnotesize$2\d\!-\!\alpha$}
\put(308,28){\circle*{5}}
\put(338,33){\footnotesize$3\d\!-\!\alpha$}
\put(352,28){\circle*{5}}
\put(382,33){\footnotesize$4\d\!-\!\alpha$}
\put(396,28){\circle*{5}}
\put(420,28){\ldots}
\put(0,6){\ldots}
\put(24,11){\footnotesize$-\!4\d\!-\!2\alpha$}
\put(44,6){\circle{5}}
\put(112,11){\footnotesize$-\!2\d\!-\!2\alpha$}
\put(132,6){\circle{5}}
\put(212,11){\footnotesize$-\!2\alpha$}
\put(220,6){\circle{5}}
\put(294,11){\footnotesize$2\d\!-\!2\alpha$}
\put(308,6){\circle{5}}
\put(382,11){\footnotesize$4\d\!-\!2\alpha$}
\put(396,6){\circle{5}}
\put(420,6){\ldots}
\end{picture}
\end{center}

\vskip -5pt
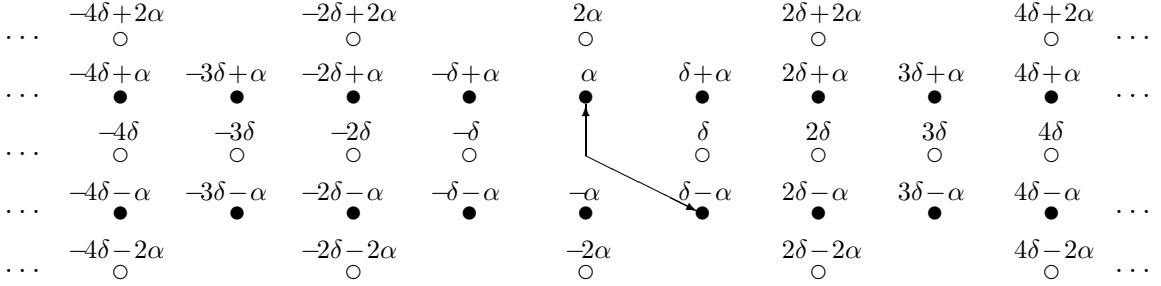
\centerline{\footnotesize Fig. 2a. The total root system $\D$
of $C(2)^{(2)}$}

\vskip 10pt
\begin{center}
\begin{picture}(450,120)
\put(0,72){\ldots}
\put(24,78){\footnotesize$-\!4\d\!+\!\alpha$}
\put(44,72){\circle*{5}}
\put(68,78){\footnotesize$-\!3\d\!+\!\alpha$}
\put(88,72){\circle*{5}}
\put(112,78){\footnotesize$-\!2\d\!+\!\alpha$}
\put(132,72){\circle*{5}}
\put(161,78){\footnotesize$-\!\d\!+\!\alpha$}
\put(176,72){\circle*{5}}
\put(218,78){\footnotesize$\alpha$}
\put(220,72){\circle*{5}}
\put(255,78){\footnotesize$\d\!+\!\alpha$}
\put(264,72){\circle*{5}}
\put(294,78){\footnotesize$2\d\!+\!\alpha$}
\put(308,72){\circle*{5}}
\put(338,78){\footnotesize$3\d\!+\!\alpha$}
\put(352,72){\circle*{5}}
\put(382,78){\footnotesize$4\d\!+\!\alpha$}
\put(396,72){\circle*{5}}
\put(420,72){\ldots}
\put(0,50){\ldots}
\put(35,56){\footnotesize$-\!4\d$}
\put(44,50){\circle{5}}
\put(79,56){\footnotesize$-\!3\d$}
\put(88,50){\circle{5}}
\put(123,56){\footnotesize$-\!2\d$}
\put(132,50){\circle{5}}
\put(169,56){\footnotesize$-\!\d$}
\put(176,50){\circle{5}}
\put(220,50){\vector(0,1){19}}
\put(220,50){\vector(2,-1){42}}
\put(262,56){\footnotesize$\d$}
\put(264,50){\circle{5}}
\put(303,56){\footnotesize$2\d$}
\put(308,50){\circle{5}}
\put(347,56){\footnotesize$3\d$}
\put(352,50){\circle{5}}
\put(391,56){\footnotesize$4\d$}
\put(396,50){\circle{5}}
\put(4200,50){\ldots}
\put(0,28){\ldots}
\put(24,33){\footnotesize$-\!4\d\!-\!\alpha$}
\put(44,28){\circle*{5}}
\put(68,33){\footnotesize$-\!3\d\!-\!\alpha$}
\put(88,28){\circle*{5}}
\put(112,33){\footnotesize$-\!2\d\!-\!\alpha$}
\put(132,28){\circle*{5}}
\put(161,33){\footnotesize$-\!\d\!-\!\alpha$}
\put(176,28){\circle*{5}}
\put(213,33){\footnotesize$-\!\alpha$}
\put(220,28){\circle*{5}}
\put(255,33){\footnotesize$\d\!-\!\alpha$}
\put(264,28){\circle*{5}}
\put(294,33){\footnotesize$2\d\!-\!\alpha$}
\put(308,28){\circle*{5}}
\put(338,33){\footnotesize$3\d\!-\!\alpha$}
\put(352,28){\circle*{5}}
\put(382,33){\footnotesize$4\d\!-\!\alpha$}
\put(396,28){\circle*{5}}
\put(420,28){\ldots}
\end{picture}
\end{center}

\vskip -10pt
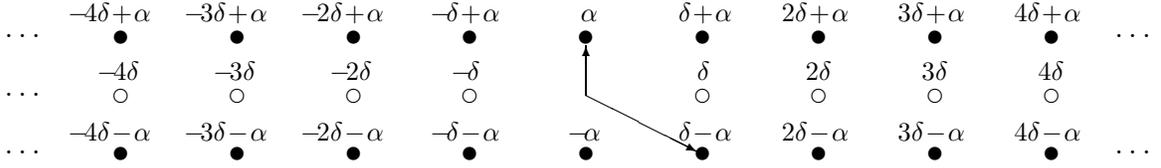
\centerline{\footnotesize Fig. 2b. The reduced root
system $\uD$ of $C(2)^{(2)}$}.

\setcounter{equation}{0}
\section{Defining relations of $U_q(A_1^{(1)})$ and
$U_q(C(2)^{(2)})$}
The quantum (q-deformed) affine (super)algebras $U_q(A_1^{(1)})$
and $U_q(C(2)^{(2)})$ are generated by the Chevalley elements
$k_{\rm d}^{\pm1}:=q^{\pm h_{\rm d}}$,
$k_{\a}^{\pm1}:=q^{\pm h_{\a}}$, $k_{\d-\a}^{\pm1}:=q^{\pm h_{\d-\a}}$,
$e_{\pm\a}$, $e_{\pm(\d-a)}$ with the defining relations
\bneqn
k_{\g}^{}k_{\g}^{-1}&\!\!=\!\!&k_{\g}^{-1}k_{\g}^{}=1~,\qquad\quad\;\;
[k_{\g}^{\pm 1},k_{\g^{\prime}}^{\pm 1}]=0~,\qquad
\lb{DR1}
\\[7pt]
k_{\g}^{}e_{\pm\a}^{}k_{\g}^{-1}&\!\!=\!\!&
q^{\pm(\g,\a)}e_{\pm\a}^{}~,\qquad\;
k_{\g}^{}e_{\pm(\d-\a)}^{}k_{\g}^{-1}=
q^{\pm(\g,\d-\a)}e_{\pm(\d-\a)}^{}~,
\lb{DR2}
\\[7pt]
[e_{\a}^{},e_{-\a}^{}]&\!\!=\!\!&[h_{\a}]_q~,\qquad\qquad\;\;
[e_{\d-\a}^{},e_{-\d+\a}^{}]=[h_{\d-\a}]_q~,
\lb{DR3}
\\[7pt]
[e_{\a}^{},e_{-\d+\a}^{}]&\!\!=\!\!&0~,\qquad\qquad\qquad\quad
[e_{-\a}^{},e_{\d-\a}^{}]=0~,\;\phantom{possible, poss}
\lb{DR4}
\edeqn
\bneqn
[e_{\pm\a}^{},[e_{\pm\a}^{},[e_{\pm\a}^{},
e_{\pm(\d-\a)}^{}]_q]_q]_{q}&\!\!=\!\!&0~,\quad\;\;
\lb{DR5}
\\[7pt]
[[[e_{\pm\a},e_{\pm(\d-\a)}]_q,e_{\pm(\d-\a)}]_q,
e_{\pm(\d-\a)}]_q&\!\!=\!\!&0~,\quad\;\;
\lb{DR6}
\edeqn
where ($\g\!=\!{\rm d},\a,\d-\a$), $({\rm d},\a)\!=\!0$,
$({\rm d},\d)\!=\!1$, and
$[h_\b]_q\!:=\!(k_\b\!-\!k_\b^{-\!1})/(q\!-\!q^{-\!1})$.
The brackets $[\cdot,\cdot]$ and $[\cdot,\cdot]_{q}$ are the super-,
and q-super-commutators:
\bn
\begin{array}{rcl}
[e_{\b}^{},e_{\b'}^{}]&\!\!=\!\!&e_{\b}^{}e_{\b'}^{}-
(-1)^{\vth(\b)\vth(\b')}e_{\b'}^{}e_{\b}^{}~,
\\[7pt]
[e_{\b}^{},e_{\b'}^{}]_{q}&\!\!=\!\!&e_{\b}^{}e_{\b'}^{}-
(-1)^{\vth(\b)\vth(\b')}q^{(\b,\b')}e_{\b'}^{}e_{\b}^{}~.
\lb{DR7}
\end{array}
\ed
Here the symbol $\vth(\cdot)$ means the parity function: $\vth(\b)=0$
for any even root $\b$, and $\vth(\b)=1$ for any odd root $\b$.

\nin
{\it Remark}. The left sides of the relations (\ref{DR5}) and (\ref{DR6})
are invariant with respect to the replacement of $q$ by $q^{-1}$.
Indeed, if we remove the q-brackets we see that the left sides of
(\ref{DR5}) and (\ref{DR6}) contain the symmetric functions of $q$
and $q^{-1}$. This property permits to write the q-commutators in
(\ref{DR5}) and (\ref{DR6}) in the inverse order, i.e.
\bneqn
[[[e_{\pm(\d-\a)}^{},e_{\pm\a}^{}]_{q},e_{\pm\a}^{}]_{q},
e_{\pm\a}^{}]_q&\!\!=\!\!&0\ ,
\lb{DR8}
\\[7pt]
[e_{\pm(\d-\a)},[e_{\pm(\d-\a)},[e_{\pm(\d-\a)},
e_{\pm\a}]_q]_q]_q&\!\!=\!\!&0~.
\lb{DR9}
\edeqn
Now we prove the useful proposition.
\bnpr
(i) In the quantum (super)algebras $U_q(A_1^{(1)})$ and
$U_q(C(2)^{(2)})$ the following relations
\bneqn
[[e_{\pm\a}^{},[e_{\pm\a}^{},e_{\pm(\d-\a)}^{}]_q]_q,
[[e_{\pm\a}^{},[e_{\pm\a}^{},e_{\pm(\d-\a)}^{}]_q]_q,
[e_{\pm\a}^{},e_{\pm(\d-\a)}^{}]_q]_q]_q&\!\!=\!\!&0~,
\lb{DR12}
\\[7pt]
[[e_{\pm(\d-\a)}^{},[e_{\pm(\d-\a)}^{},e_{\pm\a}^{}]_q]_q,
[[e_{\pm(\d-\a)}^{},[e_{\pm(\d-\a)}^{},e_{\pm\a}^{}]_q]_q,
[e_{\pm(\d-\a)}^{},e_{\pm\a}^{}]_q]_q]_q\!\!&=\!\!&0~.
\qquad\quad{}
\lb{DR13}
\edeqn
are fulfilled.

\nin
(ii) On the contrary, if the relations (\ref{DR1})--(\ref{DR4}) and
(\ref{DR12}), (\ref{DR13}) are satisfied then also  the relations
(\ref{DR5}), (\ref{DR6}) are valid.
\lb{PDR1}
\edpr
Thus, the proposition says that under the conditions
(\ref{DR1})--(\ref{DR4}) the relations (\ref{DR5}), (\ref{DR6}) and
(\ref{DR12}), (\ref{DR13}) are equivalent.

\nin
{\it Proof}. Let us assume that the relations (\ref{DR1})--(\ref{DR6})
are fulfilled.
We take the relations (\ref{DR6}) and apply to them the
corresponding q-commutator with fourth power of $e_{\pm\a}$,
i.e. $[e_{\pm\a},[e_{\pm\a},[e_{\pm\a},
[e_{\pm\a},a_{\pm}]_q\dots]_q\!=\!0$, where $a_{\pm}$ is the
left side of (\ref{DR6}). After tedious calculation we arrive to the
relations (\ref{DR12}). The relations (\ref{DR13}) are proved in
similar way. Namely, the relations $[e_{\pm(\d-\a)},
[e_{\pm(\d-\a)},[e_{\pm(\d-\a)},[e_{\pm(\d-\a)},b_{\pm}]_q\dots]_q\!=\!0$,
where $b_{\pm}$ is the right-side of (\ref{DR5}) in the form (\ref{DR8}),
results in (\ref{DR13}). On the contrary, if the relations
(\ref{DR1})--(\ref{DR4}) and (\ref{DR12}), (\ref{DR13}) are realized
then from the relations
$[e_{\mp\a},[e_{\mp\a},[e_{\mp\a},[e_{\mp\a},a_{\mp}']\dots]\!=\!0$
and $[e_{\mp(\d+\a)},[e_{\mp(\d+\a)},[e_{\mp(\d+\a)},
[e_{\mp(\d+\a)},b_{\mp}']\dots]\!=\!0$,
where  correspondingly $a_{\mp}'$ and $b_{\mp}'$ are the left sides of
the relations (\ref{DR12}) and (\ref{DR13}), follow the relations
(\ref{DR5}) and (\ref{DR6}). $\Box$

The standard Hopf structure of the quantum (super)algebras
$U_q(A_1^{(1)})$ and $U_q(C(2)^{(2)})$ is given by the following
formulas for the comultiplication $\D_q$ and antipode $S_q$:
\bn
\begin{array}{rcccl}
\D_{q}(k_\g^{\pm1})&\!\!=\!\!&k_\g^{\pm1}\ot k_\g^{\pm1}~,
\qquad\qquad\qquad S_{q}(k_\g^{\pm1})&\!\!=\!\!&k_\g^{\mp1}~,
\\[7pt]
\D_{q}(e_{\b}^{})&\!\!=\!\!&e_{\b}^{}\ot 1+ k_{\b}^{-1}\ot e_{\b}^{}~,
\qquad\quad S_{q}(e_{\b}^{})&\!\!=\!\!&-k_{\b}e_{\b}^{}~,
\\[7pt]
\D_{q}(e_{-\b}^{})&\!\!=\!\!&e_{-\b}^{}\ot k_{\b}+1 \ot e_{-\b}^{}~,
\qquad S_{q}(e_{-\b}^{})&\!\!=\!\!&-e_{-\b}^{}k_{\b}^{-1}~,
\lb{DR15}
\end{array}
\ed
where $\b=\a,\,\d-\a$; $\g={\rm d},\,\b$.

It is not hard to verify by direct calculations for the defining
relations (\ref{DR1})--(\ref{DR6}) that the quantum affine
(super)algebras $U_q(A_1^{(1)})$ and $U_q(C(2)^{(2)})$ have
the following simple involute (anti)automorphisms.

\nin
{\it (i) The non-graded antilinear antiinvolution or conjugation
"$^{*}$"}:
\bn
\begin{array}{rcccl}
(q^{\pm1})^{*}\!\!&=\!\!&q^{\mp1}~, \qquad\quad
(k_{\g}^{\pm 1})^{*}\!\!&=\!\!&k_{\g}^{\mp 1}~,\qquad
\\[7pt]
e_{\b}^{*}\!\!&=\!\!&e_{-\b}~,\qquad\qquad
e_{-\b}^{*}\!\!&=\!\!&e_{\b}
\lb{DR16}
\end{array}
\ed
($(xy)^*=y^*x^*$ for $\forall\;\, x,y\in U_{q}(g)$).

\nin
{\it (ii) The graded antilinear antiinvolution or graded conjugation
"$^{\ddagger}$"}:
\bn
\begin{array}{rcccl}
(q^{\pm1})^{\ddagger}\!\!&=\!\!&q^{\mp1}~, \qquad\qquad
(k_{\g}^{\pm 1})^{\ddagger}\!\!&=\!\!&k_{\g}^{\mp 1}~,\qquad
\\[7pt]
e_{\b}^{\ddagger}\!\!&=\!\!&(-1)^{\vth(\b)}e_{-\b}~,\quad\;
e_{-\b}^{\ddagger}\!\!&=\!\!&e_{\b}
\lb{DR17}
\end{array}
\ed
($(xy)^{\ddagger}=(-1)^{\deg x\,\deg y} y^{\ddagger}x^{\ddagger}$
for any homogeneous elements $x,y\in U_q(g)$).

\nin
{\it (iii) The Chevalley graded involution $\om$}:
\bn
\begin{array}{rcccl}
{}\qquad\om(q^{\pm1})&\!\!=\!\!&q^{\mp1}~,\qquad\quad
\om(k_{\g}^{\pm1})&\!\!=\!\!&k_{\g}^{\pm1}~,
\\[7pt]
\om(e_{\b})&\!\!=\!\!&-e_{-\b}~,\qquad\;
\om(e_{-\b})&\!\!=\!\!&-(-1)^{\th(\b)}e_{\b}~.
\lb{DR18}
\end{array}
\ed

\nin
{\it (iv) The Dynkin involution $\t$} which is associated with the
automorphism of the Dynkin diagrams of the (super)algebras
$A_1^{(1)}$ and $C(2)^{(2)}$:
\bn
\begin{array}{rcccl}
\t(q^{\pm1})&\!\!=\!\!&q^{\pm1}~,\qquad\quad
\t(k_{\rm d}^{\pm1})&\!\!=\!\!&k_{\rm d}^{\pm1}~,
\\[7pt]
\t(k_{\b}^{\pm1})&\!\!=\!\!&k_{\d-\b}^{\pm1}~,\qquad\;
\t(k_{-\b}^{\pm1})&\!\!=\!\!&k_{-\d+\b}^{\pm1}~,
\\[7pt]
\t(e_{\b})&\!\!=\!\!&e_{\d-\b}~,\qquad\;
\t(e_{-\b})&\!\!=\!\!&e_{-\d+\b}~.
\lb{DR19}
\end{array}
\ed
Here in (\ref{DR16})--(\ref{DR19}) $\b\!=\!\a,\d\!-\!\a$;
$\g\!=\!{\rm d},\b$.

It should be noted that the graded conjugation
"$^{\ddagger}$" and the Chevalley graded involution $\om$ are
involute (anti)automorphism of the fourth order, i.e., for example,
$(\om)^4\!=\!\id$. Note also that the Dynkin involution $\t$
commutes with all other three involutions, i.e. $\t(x^*)\!=\!(\t(x))^*$,
$\t(x^{\ddagger})\!=\!(\t(x))^{\ddagger}$ and $\om\t(x)\!=\!\t\om(x)$
for any element $x\in U_q(g)$ ($g\!=\!A_1^{(1)},C(2,0)^{(2)}$).

In the next Section we consider a q-analog of automorphisms
connected with the Weyl group of the (super)algebras $A_1^{(1)}$
and $C(2,0)^{(2)}$. This q-analog defines actions of the braid
group associated with the Weyl group.

\setcounter{equation}{0}
\section{Braid group actions}
We introduce the morphisms $T_{\a}$ and $T_{\d-\a}$ defined by
the following formulas:
\bn
\begin{array}{rcccl}
T_{\a}(q^{\pm1})&\!\!=\!\!&q^{\pm1}~,\qquad\qquad\qquad\qquad
T_{\a}(k_{\g}^{\pm1})&\!\!=\!\!&k_{\g}^{\pm1}
k_{\a}^{\mp\frac{2(\a,\g)}{(\a,\a)}}~,
\\[7pt]
T_{\a}(e_{\a})&\!\!=\!\!&-e_{-\a}k_{\a}~,\qquad\qquad\qquad
T_{\a}(e_{-\a})&\!\!=\!\!&-(-1)^{\th(\a)}k_{\a}^{-1}e_{\a}~,\qquad\;
\\[7pt]
T_{\a}(e_{\d-a})&\!\!=\!\!&\mbox{\ls$\frac{1}{a}$}
[e_{\a},[e_{\a},e_{\d-a}]_q]_q~,\qquad\qquad\quad\;\;\;&&
\\[7pt]
T_{\a}(e_{-\d+a)})&\!\!=\!\!&\mbox{\ls$\frac{(-1)^{\th(\a)}}{a}$}
[[e_{-\d+\a},e_{-\a}]_{q^{-1}},e_{-\a}]_{q^{-1}}~,&&
\lb{WG1}
\end{array}
\ed
\bn
\begin{array}{rcccl}
T_{\d-\a}(q^{\pm1})&\!\!=\!\!&q^{\pm1}~,\qquad\qquad\qquad\quad\;
T_{\d-\a}(k_{\g}^{\pm1})&\!\!=\!\!&k_{\g}^{\pm1}
k_{\d-\a}^{\mp\frac{2(\d-\a,\g)}{(\a,\a)}},
\\[7pt]
T_{\d-\a}(e_{\d-\a})&\!\!=\!\!&-e_{-\d+\a}k_{\d-a}~,\qquad\quad
T_{\d-\a}(e_{-\d+\a})&\!\!=\!\!&-(-1)^{\th(\a)}
k_{\d-a}^{-1}e_{\d-\a},
\\[7pt]
T_{\d-\a}(e_{\a})&\!\!=\!\!&\mbox{\ls$\frac{1}{a}$}
[e_{\d-\a},[e_{\d-\a},e_{\a}]_q]_q~,\qquad\qquad\quad&&
\\[7pt]
T_{\d-\a}(e_{-\a})&\!\!=\!\!&\mbox{\ls$\frac{(-1)^{\th(\a)}}{a}$}
[e_{-\a},e_{-\d+\a}]_{q^{-1}},e_{-\d+\a}]_{q^{-1}},\!&&
\lb{WG2}
\end{array}
\ed
where $\g={\rm d},\a,\d-\a$.
It is not difficult to prove by direct verification that the
morphisms $T_{\a}^{-1}$ and $T_{\d-\a}^{-1}$ given by
\bn
\begin{array}{rcccl}
T_{\a}^{-1}(q^{\pm1})&\!\!=\!\!&q^{\pm1}~,\qquad\qquad\qquad\qquad
\quad T_{\a}^{-1}(k_{\g}^{\pm1})&\!\!=\!\!&k_{\g}^{\pm1}
k_{\a}^{\mp\frac{2(\a,\g)}{(\a,\a)}},\qquad\quad
\\[7pt]
T_{\a}^{-1}(e_{\a})&\!\!=\!\!&-(-1)^{\th(\a)}k_{\a}^{-1}e_{-\a}^{}~,
\qquad\quad\; T_{\a}^{-1}(e_{-\a}^{})&\!\!=\!\!&-e_{\a}^{}k_{\a}^{}~,
\\[7pt]
T_{\a}^{-1}(e_{\d-a}^{})&\!\!=\!\!&\mbox{\ls$\frac{1}{a}$}
[[e_{\d-\a}^{},e_{\a}^{}]_q,e_{\a}^{}]_q~,\qquad\qquad\qquad\quad\;&&
\\[7pt]
T_{\a}^{-1}(e_{-\d+a)}^{})&\!\!=\!\!&\mbox{\ls$\frac{(-1)^{\th(\a)}}{a}$}
[e_{-\a}^{},[e_{-\a}^{},e_{-\d+\a}^{}]_{q^{-1}}]_{q^{-1}}~,
\quad\;\;&&
\lb{WG3}
\end{array}
\ed
\bn
\begin{array}{rcccl}
T_{\d-\a}^{-1}(q^{\pm1})&\!\!=\!\!&q^{\pm1}~,\qquad\qquad\qquad\qquad
\quad T_{\d-\a}^{-1}(k_{\g}^{\pm1})&\!\!=\!\!&k_{\g}^{\pm1}
k_{\d-\a}^{\mp\frac{2(\d-\a,\g)}{(\a,\a)}},
\\[7pt]
T_{\d-\a}^{-1}(e_{\d-\a}^{})&\!\!=\!\!&
-(-1)^{\th(\a)}k_{\d-a}^{-1}e_{-\d+\a}^{}~,\quad\;
T_{\d-\a}^{-1}(e_{-\d+\a}^{})&\!\!=\!\!&-e_{\d-\a}^{}
k_{\d-a}^{}~,\qquad\;\;
\\[7pt]
T_{\d-\a}^{-1}(e_{\a}^{})&\!\!=\!\!&\mbox{\ls$\frac{1}{a}$}
[[e_{\a}^{},e_{\d-\a}^{}]_q,e_{\d-\a}^{}]_q~,\qquad\qquad\qquad\quad&&
\\[7pt]
T_{\d-\a}^{-1}(e_{-\a}^{})&\!\!=\!\!&\mbox{\ls$\frac{(-1)^{\th(\a)}}{a}$}
[e_{-\d+\a}^{},[e_{-\d+\a}^{},e_{-\a}^{}]_{q^{-1}}]_{q^{-1}}
\quad\;\;&&
\lb{WG4}
\end{array}
\ed
are inverses to $T_{\a}$ and $T_{\d-\a}$, i.e.
\bn
T_{\a}^{}T_{\a}^{-1}=T_{\a}^{-1}T_{\a}^{}=1~,\qquad\quad
T_{\d-\a}^{}T_{\d-\a}^{-1}=T_{\d-\a}^{-1}T_{\d-\a}^{}=1~.
\lb{WG5}
\ed
Here in (\ref{WG1})--(\ref{WG4}) and in what follows we use the
notation:
\bn
a:=[(\a,\a)]_q=\mbox{\ls$\frac{q^{(\a,\a)}-q^{-(\a,\a)}}{q-q^{-1}}$}~.
\lb{DR14}
\ed
\bnpr
(i) The morphisms $T_{\a}$ and $T_{\d-\a}$ (and also $T_{\a}^{-1}$ and
$T_{\d-\a}^{-1}$) commute with the graded conjugation "$^{\ddagger}$",
i.e.
\bn
(T_{\a}(x))^{\ddagger}=T_{\a}(x^{\ddagger})~,\qquad
(T_{\d-\a}(x))^{\ddagger}=T_{\d-\a}(x^{\ddagger})
\lb{WG6}
\ed
for any element $x\in U_q(g)$.

\nin
(ii) The morphisms $T_{\a}^{\pm1}$ and $T_{\d-\a}^{\pm1}$ are also
compatible with the Chevalley graded involution $\om$, in sense that:
\bn
T_{\a}^{}\om=\om T_{\a}^{-1}~,\qquad
T_{\d-\a}^{}\om=\om T_{\d-\a}^{-1}~.
\lb{WG7}
\ed
(iii) The morphisms $T_{\a}^{\pm1}$ and $T_{\d-\a}^{\pm1}$ are connected
with each other by the Dynkin involution $\t$, in sense that:
\bn
T_{\a}^{}\t=\t T_{\d-\a}^{}~, \qquad
T_{\a}^{-1}\t=\t T_{\d-\a}^{-1}~.
\lb{WG8}
\ed
\lb{PWG1}
\edpr
This proposition can be proved by direct verification for the
Chevalley basis.
\bnpr
The morphisms $T_{\a}$ and $T_{\d-\a}$ (and also $T_{\a}^{-1}$ and
$T_{\d-\a}^{-1}$) are the automorphisms of the quantum
(super)algebras $U_q(A_1^{(1)})$ and $U_q(C(2)^{(2)})$.
\lb{PWG2}
\edpr
{\it Proof}. The proposition is proved by direct verification that
the defining relations remain valid under the actions of the given
morphisms. To this end we apply Proposition \ref{PDR1}.
Note that under the action of $T_{\a}$ the relations (\ref{DR4}) and
(\ref{DR5}) are transformed into each other, the relation
(\ref{DR6}) is transformed to (\ref{DR12}). Analogously, under
action of $T_{\d-\a}$ the relations (\ref{DR4}) and (\ref{DR6}) are
transformed into each other, the relations
(\ref{DR5}) are transformed to (\ref{DR13}). In addition, it is useful
to apply the relations (\ref{WG6}). $\Box$

It is easy to see that all these automorphisms $T_{\a}^{\pm1}$ and
$T_{\d-\a}^{\pm1}$ are not Hopf algebra automorphisms of
$U_q(A_1^{(1)})$ and $U_q(C(2)^{(2)})$, in sense that, e.g.,
$T_{\a}\ot T_{\a}\circ\D_q\neq \D_q\circ T_{\a}$.

In the case of $U_q(g)$, where $g$ is a finite-dimensional simple
Lie algebras, the automorphisms of type $T_{\a}^{\pm1}$ and
$T_{\d-\a}^{\pm1}$ are called the Lusztig automorphisms \cite{L}.

Introduce the following root vectors:
\bn
\begin{array}{rcccl}
e_{\d}\!\!&:=\!\!&[e_{\a},e_{\d-\a}]_q~,\qquad
e_{-\d}\!\!&:=\!\!&[e_{-\d+\a},e_{-\a}]_{q^{-1}},
\\[7pt]
\tl{e}_{\d}\!\!&:=\!\!&[e_{\d-\a},e_{\a}]_q~, \qquad
\tl{e}_{-\d}\!\!&:=\!\!&[e_{-\a},e_{-\d+\a}]_{q^{-1}}.
\lb{WG9}
\end{array}
\ed
It is not difficult to verify that under the actions of the automorphisms
$T_{\a}^{\pm1}$ and $T_{\d-\a}^{\pm1}$ the elements $e_{\pm\d}$ and
$\tl{e}_{\pm\d}$ are transformed as follows:
\bn
\begin{array}{rcccl}
T_{\a}^{}(\tl{e}_{\pm\d}^{})\!\!&=\!\!&
(-1)^{\th(\a)}e_{\pm\d}^{}~,\qquad\;
T_{\a}^{-1}(e_{\pm\d}^{})\!\!&=\!\!&
(-1)^{\th(\a)}\tl{e}_{\pm\d}^{}~,
\\[7pt]
T_{\d-\a}^{}(e_{\pm\d}^{})\!\!&=\!\!&
(-1)^{\th(\a)}\tl{e}_{\pm\d}^{}~,\qquad
T_{\d-\a}^{-1}(\tl{e}_{\pm\d}^{})\!\!&=\!\!&
(-1)^{\th(\a)}e_{\pm\d}^{}~.
\lb{WG10}
\end{array}
\ed
Therefore
\bn
\begin{array}{rcccl}
T_{2\d}^{}(e_{\pm\d}^{})\!\!&=\!\!&e_{\pm\d}^{}~,\qquad
T_{2\d}^{-1}(e_{\pm\d}^{})\!\!&=\!\!&e_{\pm\d}^{}~,
\\[7pt]
\tl{T}_{2\d}^{}(\tl{e}_{\pm\d}^{})\!\!&=\!\!&
\tl{e}_{\pm\d}^{}~,\qquad
\tl{T}_{2\d}^{-1}(\tl{e}_{\pm\d}^{})\!\!&=\!\!&
\tl{e}_{\pm\d}^{}~,
\lb{WG11}
\end{array}
\ed
where the elements $T_{2\d}$ and $\tl{T}_{2\d}$ called
the translation operators are given by
\bn
\begin{array}{rcccl}
T_{2\d}^{}\!\!&=\!\!&T_{\a}^{}T_{\d-\a}^{}~,\qquad
T_{2\d}^{-1}\!\!&=\!\!&T_{\d-\a}^{-1}T_{a}^{-1}~,
\\[7pt]
\tl{T}_{2\d}^{}\!\!&=\!\!&T_{\d-\a}^{}T_{\a}^{}~,\qquad
\tl{T}_{2\d}^{-1}\!\!&=\!\!&T_{\a}^{-1}T_{\d-\a}^{-1}~.
\lb{WG12}
\end{array}
\ed
\bnpr
The automorphisms $T_{\d}^{}\!:=\!T_{\a}^{}\t$,
$T_{\d}^{-1}\!:=\!\t T_{\a}^{-1}$ and
$\tl{T}_{\d}^{}\!:=\!T_{\d-\a}^{}\t$,
$\tl{T}_{\d}^{-1}\!:=\!\t T_{\d-\a}^{-1}$ are the square roots of
the automorphisms of $T_{2\d}^{\pm1}$ and $\tl{T}_{2\d}^{\pm1}$
correspondingly, i.e.
\bn
\begin{array}{rcccl}
T_{\d}^2\!\!&=\!\!&T_{2\d}^{}~,\qquad
T_{\d}^{-2}\!\!&=\!\!&\tl{T}_{2\d}^{-1}~,
\\[7pt]
\tl{T}_{\d}^2\!\!&=\!\!&\tl{T}_{2\d}~,\qquad
\tl{T}_{\d}^{-2}\!\!&=\!\!&\tl{T}_{2\d}^{-1}~.
\lb{WG13}
\end{array}
\ed
Moreover
\bn
\begin{array}{rcccl}
T_{\d}^{}(e_{\d-\a}^{})\!\!&=\!\!&-e_{-\a}^{}k_{\a}^{}~,
\qquad\qquad\qquad
T_{\d}^{}(e_{-\d+\a}^{})\!\!&=\!\!&-(-1)^{\th(\a)}
k_{\a}^{-1}e_{\a}^{},
\\[7pt]
T_{\d}^{-1}(e_{\a}^{})\!\!&=\!\!&-(-1)^{\th(\a)}
k_{\d-\a}^{-1}e_{-\d+\a}^{}~,\quad
T_{\d}^{-1}(e_{-\a}^{})\!\!&=\!\!&-e_{\d-\a}^{}k_{\d-\a}^{}~,
\lb{WG14}
\\[7pt]
\end{array}
\ed
\bn
\begin{array}{rcccl}
\tl{T}_{\d}^{}(e_{\a}^{})\!\!&=\!\!&
-e_{-\d+\a}^{}k_{\d-\a}^{}~,\qquad\qquad\quad
\tl{T}_{\d}(e_{-\a})\!\!&=\!\!&-(-1)^{\th(\a)}
k_{\d-\a}^{-1}e_{\d-\a}^{},\;
\\[7pt]
\tl{T}_{\d}^{-1}(e_{\d-\a})\!\!&=\!\!&-(-1)^{\th(\a)}
k_{\a}^{-1}e_{-\a}^{}~,\quad
\tl{T}_{\d}^{-1}(e_{-\d+\a}^{})\!\!&=\!\!&-e_{\a}^{}k_{\a}~,
\lb{WG15}
\end{array}
\ed
and also
\bn
\begin{array}{rcccl}
T_{\d}^{}(e_{\pm\d})\!\!&=\!\!&(-1)^{\th(\a)}e_{\pm\d}~,\qquad\;
T_{\d}^{-1}(e_{\pm\d})\!\!&=\!\!&(-1)^{\th(\a)}e_{\pm\d}~,\qquad\,
\\[7pt]
\tl{T}_{\d}^{}(\tl{e}_{\pm\d})\!\!&=\!\!&
(-1)^{\th(\a)}\tl{e}_{\pm\d}~,\qquad\;\;
\tl{T}_{\d}^{-1}(\tl{e}_{\pm\d})\!\!&=\!\!&
(-1)^{\th(\a)}\tl{e}_{\pm\d}~.
\lb{WG16}
\end{array}
\ed
\lb{PWG3}
\edpr
{\it Proof}. From ({\ref{WG8}) we have that
$T_{\d-\a}\!=\!\t T_{\a}\t$ and therefore
$T_{2d}=T_{\a}T_{\d-\a}\!=\!T_{\a}\t T_{\a}\t=T_{\d}^2$.
Analogously $T_{\a}\!=\!\t T_{\d-\a}\t$ and therefore $\tl{T}_{2d}\!
=\!T_{\d-\a}T_{\a}\!=\!T_{\d-\a}\t T_{\d-\a}\t=\tl{T}_{\d}^2$.
The formulas ({\ref{WG14})--(\ref{WG16}) are trivial.

In the next section we construct the Cartan-Weyl basis and describe
its properties in detail.

\setcounter{equation}{0}
\section{Cartan-Weyl basis for $U_q(A_1^{(1)})$ and
$U_q(C(2)^{(2)})$}
A general scheme for construction of a Cartan-Weyl basis for
quantized Lie algebras and superalgebras was proposed in
Ref. \cite{T2}. The scheme was applied in detail at first for
quantized finite-dimensional Lie (super)algebras \cite{KT1} and
then to quantized non-twisted affine algebras \cite{TK}.

This procedure is bases on a notion of ``normal ordering'' for the
reduced positive root system. For affine Lie (super)algebras this
notation was formulated in \cite{T1} (see also \cite{T2},
\cite{KT2}--\cite{KT5}).
In our case the reduced positive system has only two normal
orderings:
\bn
\mbox{\sms$\a,\d\!+\!\a,2\d\!+\!\a,\ldots,\infty\d\!+\!\a,
\d,2\d,3\d,\ld,\infty\d,\infty\d\!-\!\a,\ld,
3\d\!-\!\a,2\d\!-\!\a,\d\!-\!\a$},
\lb{CW1}
\ed
\bn
\mbox{\sms$\d\!-\!\a,2\d\!-\!\a,3\d\!-\!\a,\ld,
\infty\d\!-\!\a,\d,2\d,3\d,\ldots,\infty\d,\infty\d\!+\!\a,
\ld,2\d\!+\!\a,\d\!+\!\a,\a$}.
\lb{CW2}
\ed

\nin
The first normal ordering (\ref{CW1}) corresponds to ``clockwise''
ordering for positive roots in Fig. 1, 2b if we start from root $\a$
to root $\d\!-\!\a$. The inverse normal ordering (\ref{CW2})
corresponds to ``anticlockwise'' ordering for the positive roots
when we move from $\d\!-\!\a$ to $\a$.

In accordance with the normal ordering (\ref{CW1}) we set
\bneqn
e_{\d}^{}\!\!\!&:=\!\!\!&[e_{\a}^{},e_{\d-\a}^{}]_q~,\qquad\qquad
\qquad e_{-\d}^{}:=[e_{-\d+\a}^{},e_{-\a}^{}]_{q^{-1}}~,
\lb{CW3}
\\[7pt]
e_{n\d+\a}^{}\!\!\!&:=\!\!\!&\mbox{\ls$\frac{1}{a}$}
[e_{(n-1)\d+\a}^{},e_{\d}^{}]~,\qquad\;\,
e_{-n\d-\a}^{}:=\mbox{\ls$\frac{1}{a}$}
[e_{-\d}^{},e_{-(n-1)\d-\a}^{}]~,
\lb{CW4}
\\[7pt]
e_{(n+1)\d-\a}^{}\!\!\!&:=\!\!\!&\mbox{\ls$\frac{1}{a}$}
[e_{\d}^{},e_{n\d-\a}^{}]~,\qquad\;\,
e_{-(n+1)\d+\a}^{}:=\mbox{\ls$\frac{1}{a}$}
[e_{-n\d+\a}^{},e_{-\d}^{}]~,
\lb{CW5}
\\[7pt]
e_{n\d}'\!\!\!&:=\!\!\!&[e_{\a}^{},e_{n\d-\a}^{}]_{q}~,\qquad
\qquad\quad e_{-n\d}':=[e_{-n\d+\a}^{},e_{-\a}^{}]_{q^{-1}}~,
\lb{CW6}
\edeqn
where $n=1,2,\ldots$, and $a$ is given by the formula (\ref{DR14}).
Analogously for the inverse normal ordering (\ref{CW2}) we set
\bneqn
\tl{e}_{\d}^{}\!\!\!&:=\!\!\!&
[e_{\d-\a}^{},e_{\a}^{}]_q~,\qquad\qquad\qquad
\tl{e}_{-\d}^{}:=[e_{-\a}^{},e_{-\d+\a}^{}]_{q^{-1}},
\lb{CW7}
\\[7pt]
\tl{e}_{n\d+\a}^{}\!\!\!&:=\!\!\!&\mbox{\ls$\frac{1}{a}$}
[\tl{e}_{\d}^{},\tl{e}_{(n-1)\d+\a}^{}]~,\qquad\;\,
\tl{e}_{-n\d-\a}^{}:=\mbox{\ls$\frac{1}{a}$}
[\tl{e}_{-(n-1)\d-\a}^{},\tl{e}_{-\d}^{}]~,
\lb{CW8}
\\[7pt]
\tl{e}_{(n+1)\d-\a}^{}\!\!\!&:=\!\!\!&\mbox{\ls$\frac{1}{a}$}
[\tl{e}_{n\d-\a}^{},\tl{e}_{\d}^{}]~,\qquad\;\,
\tl{e}_{-(n+1)\d+\a}^{}:=\mbox{\ls$\frac{1}{a}$}
[\tl{e}_{-\d}^{},\tl{e}_{-n\d+\a}^{}]~,
\lb{CW9}
\\[7pt]
\tl{e}_{n\d}'\!\!\!&:=\!\!\!&
[e_{\d-\a}^{},\tl{e}_{(n-1)\d+\a}^{}]_{q}~,\qquad\;\;
\tl{e}_{-n\d}':=[e_{-\d+\a}^{},\tl{e}_{-(n-1)\d-\a}^{}]_{q^{-1}},
\lb{CW10}
\edeqn
where $n=1,2,\ldots$.
Thus, we have two systems of the Cartan-Weyl generators: 'direct'
and 'inverse'. Each such system together with the Cartan generators
$k_{\a}^{\pm1}$, $k_{\d-\a}^{\pm1}$, $e_{\pm\a}^{}$,
$e_{\pm(\d-\a)}^{}$ are called the q-analog of
the Cartan-Weyl basis (or simply the Cartan-Weyl basis) for
the quantum (super)algebras $U_q(A_1^{(1)})$ and $U_q(C(2)^{(2)})$}.

Now we consider some properties of these bases.
First of all, the explicit construction of the Cartan-Weyl generators
(\ref{CW3})--(\ref{CW6}) (or (\ref{CW7})--(\ref{CW10})) permits easy
to find their properties with respect to the (anti)involutions
(\ref{DR16})--(\ref{DR18}). For example, it is evident that
\bn
e_{\pm\g}^{\,*}=e_{\mp\g}~,\qquad \forall \g\,\in \uDp~.
\lb{CW11}
\ed
and also
\bn
\begin{array}{rcccl}
e_{n\d+\a}^{\,\ddagger}\!\!\!&=\!\!\!&
(\!-1)^{(n+1)\th(\a)}e_{-n\d-\a}~,\quad
e_{-n\d-\a}^{\,\ddagger}\!\!\!&=\!\!\!&
(\!-1)^{n\th(\a)}e_{n\d+\a}~,
\\[7pt]
e_{n\d-\a}^{\,\ddagger}\!\!\!&=\!\!\!&
(\!-1)^{n\th(\a)}e_{-n\d+\a}~,\quad\quad\;\;
e_{-n\d+\a}^{\,\ddagger}\!\!\!&=\!\!\!&
(\!-1)^{(n-1)\th(\a)}e_{n\d-\a}~,
\\[7pt]
e_{n\d}^{\,\ddagger}\!\!\!&=\!\!\!&
(\!-1)^{n\th(\a)}e_{n\d}~,\qquad\qquad\quad
e_{-n\d}^{\,\ddagger}\!\!\!&=\!\!\!&
(\!-1)^{n\th(\a)}e_{-n\d}~.
\lb{CW11'}
\end{array}
\ed
Further, it is easy to see that the 'direct' and 'inverse' Cartan-Weyl
generators ({\ref{CW3})--(\ref{CW6}) and ({\ref{CW7})--(\ref{CW10})
have very simple connection by the Dynkin involution $\t$:
\bn
\begin{array}{rccccl}
\t(e_{n\d+\a})\!\!\!&=\!\!\!&\tl{e}_{(n+1)\d-\a}~,\quad\;
\t(\tl{e}_{n\d+\a})\!\!\!&=\!\!\!&e_{(n+1)\d-\a}\quad &(n\in\ZZ\,)~,\quad
\\[7pt]
\t(e_{n\d-\a})\!\!\!&=\!\!\!&\tl{e}_{(n-1)\d+\a}~,\quad\;
\t(\tl{e}_{n\d-\a})\!\!&=\!\!&e_{(n-1)\d+\a}\quad &(n\in\ZZ\,)~,
\\[7pt]
\t(e_{n\d})\!\!\!&=\!\!\!&\tl{e}_{n\d}~,\qquad\qquad\quad
\t(\tl{e}_{n\d})\!\!\!&=\!\!\!&e_{n\d}\qquad\qquad &(n\neq0)~.
\lb{CW12}
\end{array}
\ed
The transformation properties with respect to the automorphisms
$T_{\a}$ and $T_{\d-\a}$ can be not hard obtained with the help of
(\ref{WG1}), (\ref{WG2}), (\ref{WG10}), and they have the form
\bn
\begin{array}{rcccl}
T_{\a}(\!\tl{e}_{n\d+\a}\!)\!\!\!&=\!\!\!&
(\!-1)^{(n+1)\th(\a)}e_{n\d-\a}~,\quad
T_{\a}(\!\tl{e}_{-n\d-\a}\!)\!\!\!&=\!\!\!&
(\!-1)^{n\th(\a)}e_{-n\d+\a}~,
\\[7pt]
T_{\a}(\!\tl{e}_{n\d-\a}\!)\!\!\!&=\!\!\!&
(\!-1)^{(n-1)\th(\a)}e_{n\d+\a}~,\quad
T_{\a}(\!\tl{e}_{-n\d+\a}\!)\!\!\!&=\!\!\!&
(\!-1)^{n\th(\a)}e_{-n\d-\a}~,
\\[7pt]
T_{\a}(\!\tl{e}_{n\d}\!)\!\!\!&=\!\!\!&
(\!-1)^{n\th(\a)}e_{n\d}~,\qquad\qquad\;
T_{\a}(\tl{e}_{-n\d})\!\!\!&=\!\!\!&
(\!-1)^{n\th(\a)}e_{-n\d}~,
\lb{CW13}
\end{array}
\ed
where $n>0$, and
\bn
\begin{array}{rcccl}
T_{\d-\a}(\!e_{k\d+\a}\!)\!\!\!&=\!\!\!&
(\!-1)^{k\th(\a)}\tl{e}_{(k+2)\d-\a}~,\quad
T_{\d-\a}(\!e_{-k\d-\a}\!)\!\!\!&=\!\!\!&
(\!-1)^{(k+1)\th(\a)}\tl{e}_{-(k+2)\d+\a}~,
\\[7pt]
T_{\d-\a}(\!e_{l\d-\a}\!)\!\!\!&=\!\!\!&
(\!-1)^{l\th(\a)}\tl{e}_{(l-2)\d+\a}~,\quad\;
T_{\d-\a}(\!e_{-l\d+\a}\!)\!\!\!&=\!\!\!&
(\!-1)^{(l-1)\th(\a)}\tl{e}_{(-l+2)\d-\a}~,
\\[7pt]
T_{\d-\a}(\!e_{m\d}\!)\!\!\!&=\!\!\!&
(\!-1)^{m\th(\a)}\tl{e}_{m\d}~,\qquad\quad\;
T_{\d-\a}(e_{-m\d})\!\!\!&=\!\!\!&
(\!-1)^{m\th(\a)}\tl{e}_{-m\d}~,
\lb{CW14}
\end{array}
\ed
for $k\ge0$, $l>1$, $m>0$. As corollary of the formulas
(\ref{CW12})--(\ref{CW14}) we easy find the actions of
the translation operators $T_\d$ and $\tl{T}_\d$:
\bn
\begin{array}{rcccl}
{}\quad T_{\d}(\!e_{k\d+\a}\!)\!\!\!&=\!\!\!&
(\!-1)^{k\th(\a)}e_{(k+1)\d+\a}~,\qquad
T_{\d}(\!e_{-k\d-\a}\!)\!\!\!&=\!\!\!&
(\!-1)^{(k+1)\th(\a)}e_{-(k+1)\d-\a}~,
\\[7pt]
T_{\d}(\!e_{l\d-\a}\!)\!\!\!&=\!\!\!&
(\!-1)^{l\th(\a)}e_{(l-1)\d-\a}~,\qquad\;\;
T_{\d}(\!e_{-l\d+\a}\!)\!\!\!&=\!\!\!&
(\!-1)^{(l-1)\th(\a)}e_{(-l+1)\d+\a}~,
\\[7pt]
T_{\d}(\!e_{m\d}\!)\!\!\!&=\!\!\!&
(\!-1)^{m\th(\a)}e_{m\d}~,\qquad\qquad\;
T_{\d}(e_{-m\d})\!\!\!&=\!\!\!&
(\!-1)^{m\th(\a)}e_{-m\d}
\lb{CW15}
\end{array}
\ed
for $k\ge0$, $l>1$, $m>0$, and
\bn
\begin{array}{rcccl}
\tl{T}_{\d}(\!\tl{e}_{n\d+\a}\!)\!\!\!&=\!\!\!&
(-1)^{(n-1)\th(\a)}\tl{e}_{(n-1)\d+\a}~,\quad\;
\tl{T}_{\d}(\!\tl{e}_{-n\d-\a}\!)\!\!\!&=\!\!\!&
(-1)^{n\th(\a)}\tl{e}_{(-n+1)\d+\a}~,
\\[7pt]
\tl{T}_{\d}(\!\tl{e}_{n\d-\a}\!)\!\!\!&=\!\!\!&
(-1)^{(n+1)\th(\a)}\tl{e}_{(n+1)\d-\a}~,\quad\;
T_{\d}(\!\tl{e}_{-n\d+\a}\!)\!\!\!&=\!\!\!&
(-1)^{n\th(\a)}\tl{e}_{-(n+1)\d+\a}~,
\\[7pt]
\tl{T}_{\d}(\!\tl{e}_{n\d}\!)\!\!\!&=\!\!\!&
(-1)^{n\th(\a)}\tl{e}_{n\d}~,\qquad\qquad\qquad
\tl{T}_{\d}(\tl{e}_{-n\d})\!\!\!&=\!\!\!&
(-1)^{n\th(\a)}\tl{e}_{-n\d}~,
\lb{CW16}
\end{array}
\ed
where $n>0$. (Also see (\ref{WG14}) and (\ref{WG16})).
Using the formulas (\ref{CW15})--(\ref{CW16}) we can easy find
the actions for the inverse translation operators $T_{\d}^{-1}$,
$\tl{T}_{\d}^{-1}$ and $T_{2\d}^{-1}$, $\tl{T}_{2\d}^{-1}$.
These actions are not written here. From the relations
(\ref{CW15})--(\ref{CW16}) it is clear that the operators
$T_{\d}^{\pm1}$ and $\tl{T}_{\d}^{\pm1}$ can be used for
construction of the Cartan-Weyl generators (\ref{CW3})--(\ref{CW6})
starting from the Chevalley basis. In the case of the quantum
untwisted affine algebras the similar procedure was applied in
the paper \cite{Be1}.
\bnpr
The root vectors (\ref{CW3})--(\ref{CW6}) satisfy the following
permutation relations:
\bn
\begin{array}{rcccl}
k_{\rm d}^{}e_{n\d\pm\a}^{}k_{\rm d}^{-1}\!\!&=\!\!&
q^{n({\rm d},\d)}e_{n\d\pm\a}^{},\qquad
k_{\rm d}^{}e_{n\d}'k_{\rm d}^{-1}\!\!&=\!\!&
q^{n({\rm d},\d)}e_{n\d}',
\\[7pt]
k_{\g}^{}e_{n\d\pm\a}^{}k_{\g}^{}\!\!\!&=\!\!\!&
q^{\pm(\g,\a)}e_{n\d\pm\a}^{},\qquad
k_{\g}^{}e_{n\d}'k_{\g}^{-1}\!\!&=\!\!&e_{n\d}
\lb{CW17}
\end{array}
\ed
for any $n\,\in\ZZ$ and any $\g\,\in\uDp$, and also
\bneqn
[e_{n\d+\a},e_{-n\d-\a}]&\!\!=\!\!&(-1)^{n\th(\a)}
\mbox{\ls$\frac{k_{n\d+\a}^{}-k_{n\d+\a}^{-1}}{q-q^{-1}}$}
\qquad\qquad\; (n\ge0)~,
\lb{CW18}
\\[7pt]
[e_{n\d-\a},e_{-n\d+\a}]&\!\!=\!\!&(-1)^{(n-1)\th(\a)}
\mbox{\ls$\frac{k_{n\d-\a}^{}-k_{n\d-\a}^{-1}}{q-q^{-1}}$}
\qquad\quad (n>0)~;
\lb{CW19}
\\[7pt]
[e_{n\d+\a}^{},e_{(n+2m-\!1)\d+\a}^{}]_{q}&\!\!=\!\!&
(q_{\a}^2\!-\!1)\sum_{l=1}^{m-1}q_{\a}^{-l}
e_{(n+l)\d+\a}^{}e_{(n+2m-\!1\!-l)\d+\a}^{},
\lb{CW20}
\edeqn
\bn
\begin{array}{rcl}
[e_{n\d+\a}^{},e_{(n+2m)\d+\a}^{}]_{q}&\!\!=\!\!&
(q_{\a}\!-\!1)q_{\a}^{-m+1}e_{(n+m)\d+\a}^2+
\\[7pt]
&&+\,(q_{\a}^2\!-\!1)\sum\limits_{l=1}^{m-1}q_{\a}^{-l}
e_{(n+l)\d+\a}^{}e_{(n+2m-l)\d+\a}^{}\quad
\lb{CW21}
\end{array}
\ed
for any integers $n\geq0,\;m>0$;
\bn
[e_{(n+2m-1)\d-\a}^{},e_{n\d-\a}^{}]_{q}\!=\!
-(q_{\a}^2\!-\!1)\sum_{l=1}^{m-1}q_{\a}^{-l}
e_{(n+2m-\!1\!-l)\d-\a}^{}e_{(n+l)\d-\a}^{},\qquad
\lb{CW22}
\ed
\bn
\begin{array}{rcl}
[e_{(n+2m)\d-\a}^{},e_{n\d-\a}^{}]_{q}&\!\!=\!\!&
-(q_{\a}\!-\!1)q_{\a}^{-m+\!1}e_{(n+m)\d-\a}^2-
\\[7pt]
&&-\,(q_{\a}^2\!-\!1)\sum\limits_{l=1}^{m-1}q_{\a}^{-l}
e_{(n+l)\d-\a}^{}e_{(n+2m-l)\d-\a}^{} \qquad
\lb{CW23}
\end{array}
\ed
for any integers $n,\;m>0$;
\bn
\begin{array}{rcl}
[e_{-n\d+\a}^{},e_{(n+2m-\!1)\d+\a}^{}]&\!\!=\!\!&
-(-\!1)^{(n-1)\th(\a)} (q_{\a}^2\!-\!1)\,\times
\\[7pt]
&&\times\sum\limits_{l=n}^{n+m-1}q_{\a}^{-l}k_{n\d-\a}
e_{(l-n)\d+\a}^{}e_{(n+2m-\!1\!-l)\d+\a}^{}+
\\[7pt]
&&+(q_{\a}^2\!-\!1)\sum\limits_{l=1}^{n-1}(-\!1)^{l\th(\a)}
q_{\a}^{-l}k_{\d}^{l}e_{(-n+l)\d+\a}^{}e_{(n+2m-\!1\!-l)\d+\a}^{},
\end{array}
\lb{CW24}
\ed
\bn
\begin{array}{rcl}
[e_{-n\d+\a}^{},e_{(n+2m)\d+\a}^{}]&\!\!=\!\!&
-(-1)^{(n-1)\th(\a)}(q_{\a}^2\!-\!1)\,\times
\\[7pt]
&&\times\sum\limits_{l=n}^{n+m-1}\!q_{\a}^{-l}
k_{n\d-\a}e_{(-n+l)\d+\a}^{}e_{(n+2m-\!1\!-l)\d+\a}^{}+
\\[7pt]
&&+(q_{\a}^2\!-\!1)\sum\limits_{l=1}^{n-1}(-\!1)^{l\th(\a)}
q_{\a}^{-l}k_{\d}^{l}e_{(-n+l)\d+\a}^{}e_{(n+2m-\!1\!-l)\d+\a}^{}-
\\[14pt]
&&-(-1)^{(n-1)\th(\a)}(q_{\a}\!-\!1)
q_{\a}^{-m-n+\!1}k_{n\d-\a}e_{m\d+\a}^2
\end{array}
\lb{CW25}
\ed
for any integers $n,\;m\geq0$;
\bn
\begin{array}{rcl}
[e_{(n+2m-1)\d-\a}^{},e_{-n\d-\a}^{}]&\!\!=\!\!&
(-1)^{(n+1)\th(\a)} (q_{\a}^2\!-\!1)\times
\\[7pt]
&&\times\!\sum\limits_{l=n+1}^{n+m-1}\!q_{\a}^{-l}
e_{(n+2m-1-l)\d+\a}^{}e_{(l-n)\d-\a}^{}k_{n\d+\a}^{-1}-
\\[7pt]
&&-(q_{\a}^2\!-\!1)\sum\limits_{l=1}^{n-1}(-\!1)^{l\th(\a)}
q_{\a}^{-l}e_{(n+2m-1-l)\d-\a}^{}e_{(-n\!+l)\d-\a}^{}k_{\d}^{-l},
\end{array}
\lb{CW26}
\ed
\bn
\begin{array}{rcl}
[e_{(n+2m)\d-\a}^{},e_{-n\d-\a}^{}]&\!\!=\!\!&
(-1)^{(n+1)\th(\a)} (q_{\a}^2\!-\!1)\times
\\[7pt]
&&\times\sum\limits_{l=n}^{n+m-1}q_{\a}^{-l}
e_{(n+2m-l)\d-\a}^{}e_{(l-n)\d+\a}^{}k_{n\d+\a}^{-1}-
\\[7pt]
&&-(q_{\a}^2\!-\!1)\sum\limits_{l=1}^{n-1}(-\!1)^{l\th(\a)}
q_{\a}^{-l}e_{(n+2m-l)\d-\a}^{}e_{(-n+l)\d-\a}^{}k_{\d}^{-l}+
\\[14pt]
&&+(-1)^{(n-1)\th(\a)}(q_{\a}\!-\!1)
q_{\a}^{-m-n+\!1}e_{m\d-\a}^2k_{n\d+\a}^{-1}
\end{array}
\lb{CW27}
\ed
for any integers $n\geq0,\;m>0$;
\bneqn
[e_{n\d+\a}^{},e_{m\d-\a}^{}]_{q}&\!\!\!=\!\!\!&e_{(n+m)\d}'
\qquad\qquad\qquad\qquad\quad\;\;\,(n\geq 0,\;m>0)~,
\lb{CW28}
\\[9pt]
[e_{n\d+\a}^{},e_{-m\d-\a}^{}]&\!\!\!=\!\!\!&-(-1)^{(m+1)\th(\a)}
e_{(n-m)\d}'k_{m\d+\a}^{-1}\qquad\,(n>m\geq 0)~,
\lb{CW29}
\\[9pt]
[e_{-m\d+\a}^{},e_{n\d-\a}^{}]&\!\!\!=\!\!\!&-(-1)^{(m-1)\th(\a)}
k_{m\d-\a}e_{(n-m)\d}'\qquad\,(n>m>0)~,
\lb{CW30}
\\[7pt]
[e_{n\d}',e_{m\d}']&\!\!=\!\!&[e_{-n\d}',e_{-m\d}']=0
\qquad\qquad\qquad\,(n>0,\;m>0)~,
\lb{CW31}
\edeqn
\bn
[e_{n\d+\a}^{},e_{m\d}']=q_{\a}^{-m+1}ae_{(n+m)\d+\a}^{}+
(q_{\a}^2\!-\!1)\sum\limits_{l=1}^{m-1}q_{\a}^{-l}
e_{(n+l)\d+\a}^{}e_{(m-l)\d}'
\lb{CW32}
\ed
for any integers $n\geq0,\;m>0$;
\bn
[e_{m\d}',e_{n\d-\a}^{}]=q_{\a}^{-m+\!1}ae_{(n+m)\d-\a}^{}+
(q_{\a}^2\!-\!1)\sum\limits_{l=1}^{m-1}q_{\a}^{-l}
e_{(m-l)\d}'e_{(n+l)\d-\a}^{}
\lb{CW33}
\ed
for any integers $n,\;m>0$;
\bn
\begin{array}{rcl}
[e_{-n\d+\a}^{},e_{m\d}']&\!\!=\!\!&-(-1)^{(n-1)\th(\a)}
q_{\a}^{-m+1}ak_{n\d-\a}^{}e_{(m-n)\d+\a}^{}-
\\[7pt]
&&-(-1)^{(n-1)\th(\a)}(q_{\a}^2\!-\!1)k_{n\d-\a}^{}
\sum\limits_{l=n}^{m-1}q_{\a}^{-l}e_{(l-n)\d+\a}^{}e_{(m-l)\d}'+
\\[7pt]
&&+(q_{\a}^2\!-\!1)\sum\limits_{l=1}^{n-1}(-1)^{l\th(\a)}
q_{\a}^{-l}k_{\d}^{l}e_{(-n+l)\d+\a}^{}e_{(m-l)\d}'
\end{array}
\lb{CW34}
\ed
for any integers $m\ge n>0$;
\bn
\begin{array}{rcl}
[e_{-n\d+\a}^{},e_{m\d}']&\!\!=\!\!&(-1)^{m\th(\a)}
q_{\a}^{-m+1}ak_{\d}^{m}e_{(-n+m)\d+\a}^{}+
\\[7pt]
&&+(q_{\a}^2\!-\!1)\sum\limits_{l=1}^{m-1}(-1)^{l\th(\a)}
q_{\a}^{-l}k_{\d}^{l}e_{(-n+l)\d+\a}^{}e_{(m-l)\d}'
\end{array}
\lb{CW35}
\ed
for any integers $n>m>0$;
\bn
\begin{array}{rcl}
[e_{m\d}'e_{-n\d-\a}^{}]&\!\!=\!\!&-(-1)^{(n+1)\th(\a)}
q_{\a}^{-m+1}a e_{(m-n)\d-\a}^{}k_{n\d+\a}^{-1}-
\\[7pt]
&&-(-1)^{(n+1)\th(\a)}(q_{\a}^2\!-\!1)\sum\limits_{l=n+1}^{m-1}
q_{\a}^{-l}e_{(m-l)\d}'e_{(l-n)\d-\a}^{}k_{n\d+\a}^{-1}+
\\[7pt]
&&+(q_{\a}^2\!-\!1)\sum\limits_{l=1}^{n}(-1)^{l\th(\a)}
q_{\a}^{-l}e_{(m-l)\d}'e_{(-n+l)\d-\a}^{}k_{\d}^{-l}
\end{array}
\lb{CW36}
\ed
for any integers $m>n\ge0$;
\bn
\begin{array}{rcl}
[e_{m\d}'e_{-n\d-\a}^{}]&\!\!=\!\!&(-1)^{m\th(\a)}
q_{\a}^{-m+1}ae_{(-n+m)\d-\a}^{}k_{\d}^{-m}+
\\[7pt]
&&+(q_{\a}^2\!-\!1)\sum\limits_{l=1}^{m-1}\!(-1)^{l\th(\a)}
q_{\a}^{-l}e_{(m-l)\d}'e_{(-n+l)\d-\a}^{}k_{\d}^{-l}
\end{array}
\lb{CW37}
\ed
for any integers $n\ge m>0$.
\lb{PCW1}
\edpr
Here in the relations (\ref{CW20})--(\ref{CW37}) and in what follows
$q_{\a}\!:=\!(-1)^{\th(\a)}q^{(\a,\a)}$.

\nin
{\it Outline of proof}. First of all, the formulas (\ref{CW17}) are
trivial. The relations (\ref{CW18}) and (\ref{CW19}) are obtained by
application the translation operators $T_{\d}^{n}$ and $T_{\d}^{-n}$
to the relations (\ref{DR3}). Further, in terms of the generators
(\ref{CW3})--(\ref{CW6}) the relation (\ref{DR5}) means that
$[e_{\a},e_{\d+\a}]_q\!=\!0$. Applying to it the operator $T_{\d}^n$,
we obtain the relation (\ref{CW20}) for $m\!=\!1$. In the case
$m>1$ the formulas (\ref{CW20}) and (\ref{CW21}) are proved for
arbitrary $m$  by induction. If we apply the operator $T_{\d}^{-k}$
to the relations (\ref{CW20}) and (\ref{CW21}) for $n=0$, then in
the case $k<m$ we obtain the relations (\ref{CW24}) and (\ref{CW25}),
in the case $m<k<2m$ we obtain the relations which are obtained from
(\ref{CW26}) and (\ref{CW27}) by the conjugation "$^*$", and
finally for $k>2m$ we get the relations which are obtained from
(\ref{CW22}) and (\ref{CW23}) by the conjugation "$^*$". Further,
the relation (\ref{CW28}) for $n=0$ is trivial (see (\ref{CW6})).
Applying to (\ref{CW28}) with $n=0$ the operators $T_{\d}^{n}$,
we can obtain for any $n > 0$ and $m>0$ the relation (\ref{CW28})
as well as the relation ({\ref{CW29}). The relation ({\ref{CW30})
can be obtained from (\ref{CW28}) by repeated
application of
the operator $T_{\d}^{-1}$. The relations (\ref{CW32}) in
the case $n=0$ and (\ref{CW33}) in the case $n=1$ are proved by
direct verification with the help of the previous results.
Repeated application of the operators $T_{\d}^{\pm 1}$ to these relations
results in the general case $n, m>0$.
The relation (\ref{CW31}) is proved by direct verification with
the help of the relations (\ref{CW32}) and (\ref{CW33}).
At last, the relations (\ref{CW34})--(\ref{CW37}) can be obtained
from (\ref{CW32}) and (\ref{CW33}) by repeated application of the
operator $T_{\d}^{-1}$. $\Box$

The imaginary root vectors $e_{n\d}'$ do not satisfy the relations of
the type (\ref{CW18}) and therefore we introduce the new imaginary
roots vectors $e_{\pm n\d}$ by the following (Schur) relations:
\bn
e_{n\d}'=\sum_{p_{1}+2p_{2}+\ldots+np_{n}=n}\!\!\!\!\!
\mbox{\ls$\frac{\;\Bigl((-1)^{\th(\a)}
(q-q^{-1})\Bigr)^{\sum p_{i}-1}}{p_{1}!\cdots
p_{n}!}$}\;\;e_{\d}^{\;p_1}\cdots e_{n\d}^{\;p_n}.
\lb{CW38}
\ed
In terms of generating functions
\bneqn
{\cal E}'(u)&\!\!:=\!\!&(-1)^{\th(\a)}(q-q^{-1})
\sum_{n\geq 1}^{}e_{n\d}'u^{-n}~,
\lb{CW39}
\\[7pt]
{\cal E}(u)&\!\!=\!\!&(-1)^{\th(\a)}(q-q^{-1})
\sum_{n\geq 1}^{}e_{n\d}u^{-n}
\lb{CW40}
\edeqn
the relation (\ref{CW38}) may be rewritten in the form
\bn
{\cal E}'(u)=-1+\exp{\cal E}(u)
\lb{CW41}
\ed
or
\bn
{\cal E}(u)=\ln(1+{\cal E}'(u))~.
\lb{CW42}
\ed
This provides a formula inverse to (\ref{CW38})
\bn
e_{n\d}=\!\!\!\sum_{p_{1}+2p_{2}+\ldots+np_{n}=n}
\!\!\!\!\!\!\!\!\!\!\mbox{\ls$\frac{\Bigl((-1)^{\th(\a)}
(q^{-1}-q)\Bigr)^{\sum p_{i}-1}(\sum_{i=1}^{n}p_{i}-1)!}
{p_{1}!\cdots p_{n}!}$}\,(e_{\d}')^{p_1}\cdots (e_{n\d}')^{p_n}.
\lb{CW43}
\ed
The new root vectors corresponding to negative roots are
obtained by the Cartan conjugation $(^*)$:
\bn
e_{-n\d}=(e_{n\d})^{*}~.
\lb{CW44}
\ed
\bnpr
The new root vectors $e_{\pm n\d}$ satisfy the following
commutation relations:
\bneqn
[e_{n\d+\a}^{},e_{m\d}^{}]&\!\!\!=\!\!\!&(-1)^{(m-1)\th(\a)}
a(m)e_{(n+m)\d+\a}^{}\qquad\quad\;\;\,(n\ge0,\,m>0)~,
\lb{CW45}
\\[7pt]
[e_{m\d}^{},e_{n\d-\a}^{}]&\!\!\!=\!\!\!&(-1)^{(m-1)\th(\a)}
a(m)e_{(n+m)\d-\a}^{}\qquad\quad\qquad\;\;(n,\,m>0)~,
\lb{CW46}
\\[7pt]
[e_{-n\d+\a}^{},e_{m\d}^{}]&\!\!\!=\!\!\!&-(-1)^{(n+m)\th(\a)}
a(m)\,k_{n\d-\a}^{}e_{(m-n)\d+\a}^{}\quad\;\;(m\ge n>0)~,\phantom{pls}
\lb{CW47}
\\[7pt]
[e_{-n\d+\a}^{},e_{m\d}^{}]&\!\!\!=\!\!\!&(-1)^{\th(\a)}\!
a(m)\,k_{\d}^{m}e_{(-n+m)\d+\a}^{}\qquad\qquad\quad\;(n>m>0)~,
\lb{CW48}
\\[7pt]
[e_{m\d}^{},e_{-n\d-\a}^{}]&\!\!\!=\!\!\!&-(-1)^{(n+m)\th(\a)}\!
a(m)e_{(m-n)\d-\a}^{}\,k_{n\d+\a}^{-1}\qquad(m>n\ge0)~,
\lb{CW49}
\\[7pt]
[e_{m\d}^{},e_{-n\d+\a}^{}]&\!\!\!=\!\!\!&(-1)^{\th(\a)}\!
a(m)e_{(-n+m)\d-\a}^{}\,k_{\d}^{-m}\qquad\qquad\quad(n\ge m>0)~,
\lb{CW50}
\\[7pt]
[e_{n\d},e_{-m\d}]&\!\!\!=\!\!\!&\d_{nm}\,a(m)\,
\mbox{\ls$\frac{k_{\d}^{m}-k_{\d}^{-m}}{q-q^{-1}}$}
\qquad\qquad\qquad\qquad\qquad(n,\,m>0)~,
\lb{CW51}
\edeqn
where
\bn
a(m):=\mbox{\ls$\frac{q^{m(\a,\a)}-q^{-m(\a,\a)}}
{m(q-q^{-1})}$}~.
\lb{CW52}
\ed
\lb{PCW2}
\edpr
This can be proved by direct calculation, applying the relations of
Proposition (\ref{PCW1}) and the actions of the translation operators
$T_{\d}^{\pm1}$.

All the relations of Propositions (\ref{PCW1}), (\ref{PCW2}) together
with the  ones obtained from them by the conjugation describe complete
list of the permutation relations of the Cartan-Weyl bases
corresponding to the 'direct' normal ordering (\ref{CW1}).
Applying to these relations the Dynkin involution $\t$, it is easy
to obtain these results for the 'inverse' normal ordering (\ref{CW2}).

\setcounter{equation}{0}
\section{Extremal projector for $U_q(A_1^{(1)})$ and
$U_q(C(2)^{(2)})$}
A general formula for the extremal projector for quantized
contragredient Lie (super)algebras of finite growth was presented 
in Refs. \cite{T2}, \cite{KT2}, \cite{KT3}. Here we specialize this 
result to our case $U_q(g)$, where $g=A_1^{(1)},~C(2)^{(2)}$.

By definition, the extremal projector for $U_{q}(g)$ is a nonzero
element $p:=p(U_q(g))$ of the Taylor extension $T_q(g)$ of $U_q(g)$
(see Refs. \cite{T2}, \cite{KT2}, \cite{KT3}), satisfying
the equations
\bn
e_{\a}^{}\,p=p\,e_{-\a}^{}=0~,\qquad
e_{\d-\a}^{}\,p=p\,e_{-\d+\a}^{}=0~,\qquad
p^{2}=p~.
\lb{EP1}
\ed
The explicit expression of the extremal projector $p$ for
our case $U_q(g)$ can be presented as follows:
\bn
p=p_{+}^{}\,p_{0}^{}\,p_{-}^{}~,
\lb{EP2}
\ed
where the factors $p_{+}^{}$, $p_{0}^{}$ and $p_{+}^{}$ have
the following form
\bn
p_{+}^{}=\prod_{n\ge0}^{\ra}p_{n\d+\a}^{}~, \qquad
p_{0}^{}=\prod_{n\ge1}^{}p_{n\d}^{}~,\qquad
p_{-}^{}=\prod_{n\ge1}^{\la}p_{n\d-\a}^{}~.\qquad
\lb{EP3}
\ed
The elements $p_{\g}$ are given by the formula
\bneqn
p_{n\d+\a}^{}&\!\!=\!\!&\sum_{m=0}^{\infty}
\mbox{\ls$\frac{(-1)^{m}}{(m)_{\bar{q}_{\a}^{}}!}$}\,\vph_{n,m}^{\,+}
e_{-n\d-\a}^{m}e_{n\d+\a}^{m}~,
\lb{EP4}
\\[3pt]
p_{n\d}^{}&\!\!=\!\!&\sum_{m=0}^{\infty}
\mbox{\ls$\frac{(-1)^{m}}{m!}$}\,\vph_{n,m}^{\,0}
e_{-n\d}^{m}e_{n\d}^{m}~,
\lb{EP5}
\\[3pt]
p_{n\d-\a}^{}&\!\!=\!\!&\sum_{m=0}^{\infty}
\mbox{\ls$\frac{(-1)^{m}}{(m)_{\bar{q}_{\a}^{}}!}$}\,\vph_{n,m}^{\,-}
e_{-n\d+\a}^{m}e_{n\d-\a}^{m}~,
\lb{EP6}
\edeqn
where the coefficients $\vph_{m}^{+}$, $\vph_{m}^{0}$ and
$\vph_{m}^{-}$ are determined as follows:
\bneqn
\vph_{n,m}^{\,+}&\!\!=\!\!&\mbox{\ls$\frac{(-1)^{mn\th(\a)}
(q-q^{-1})^{m}q^{-m(\frac{m-1}{4}+n)(\a,\a)}}
{\prod\limits_{r=1}^{m}\Bigl(k_{n\d+\a}^{}
q^{(n+\frac{1}{2}+\frac{r}{2})(\a,\a)}
-(-1)^{(r-1)\th(\a)}k_{n\d+\a}^{-1}
q^{-(n+\frac{1}{2}+\frac{r}{2})(\a,\a)}\Bigr)}$}~,
\lb{EP7}
\\[3pt]
\vph_{n,m}^{\,0}&\!\!=\!\!&\mbox{\ls$\frac{n^{m}(q-q^{-1})^{n+m}
q^{-mn(\a,\a)}}{(q^{n(\a,\a)}-q^{-n(\a,\a)})^{m}(k_{\d}^{n}
q^{n(\a,\a)}-k_{\d}^{-n}q^{-n(\a,\a)})^{m}}$}~,
\lb{EP8}
\\[5pt]
\vph_{n,m}^{\,-}&\!\!=\!\!&\mbox{\ls$\frac{(-1)^{m(n-1)\th(\a)}
(q-q^{-1})^{m}q^{-m(\frac{m-5}{4}+n)(\a,\a)}}
{\prod\limits_{r=1}^{m}\Bigl(k_{n\d-\a}^{}
q^{(n-\frac{1}{2}+\frac{r}{2})(\a,\a)}
-(-1)^{(r-1)\th(\a)}k_{n\d-\a}^{-1}
q^{-(n-\frac{1}{2}+\frac{r}{2})(\a,\a)}\Bigr)}$}~.
\lb{EP9}
\edeqn
Here in the relations (\ref{EP4})--(\ref{EP6}) and in what follows
we use the notation $\bar{q}_{\a}\!:=\!(-1)^{\th(\a)}q^{-(\a,\a)}$,
and the symbol $(m)_{\bar{q}_{\a}}$ is defined by the formula
(\ref{UR4}).

Acting by the extremal projector $p$ on any highest weight
$U_{q}(g)$-module $M$ we obtain a space $M^{0}=pM$ of highest weight
vectors for $M$ if $pM$ has no singularities. An effective example of
application of the extremal projector for the case of the quantum
algebra $U_{q}(gl(n,\CC)$ can be found in Ref. \cite{T2}.

\setcounter{equation}{0}
\section{Universal $R$-matrix for $U_q(A_1^{(1)})$ and
$U_q(C(2)^{(2)})$}
Any quantum (super)algebra $U_{q}(g)$ is a non-cocommutative Hopf
(super)\-al\-gebra which has the intertwining operator called the
universal $R$-matrix.
By definition \cite{D}, the universal $R$-matrix for the Hopf
(super)algebra $U_q(g)$ is an invertible element $R$ of the Tylor
extension $T_q(g)\ot T_q(g)$ of $U_q(g)\ot U_q(g)$ 
(see Refs. \cite{KT3}--\cite{KT5}), satisfying the equations
\bneqn
\tl{\Delta}_q(a)&\!\!=\!\!&R\Delta_q(a)R^{-1} \qquad\quad\;\;
\forall\,\,a \in U_q(g)~,
\lb{UR1}
\\[7pt]
(\Delta_q\ot\id)R&\!\!=\!\!&R^{13}R^{23}~,\qquad
(\id\ot\Delta_q)R=R^{13}R^{12},
\lb{UR2}
\edeqn
where  $\tl{\Delta}_q$ is the opposite comultiplication:
\sloppy $\tl{\Delta}_q=\sg\Delta_q$,
$\sg(a\ot b)=(-1)^{\deg a\deg b} b\ot a$ for
all homogeneous elements $a,\,b \in U_q(g)$. In the relation
(\ref{UR2}) we use the standard notations
$R^{12}=\sum a_{i}\ot b_{i}\ot\id$,
$R^{13}=\sum a_{i}\ot\id\ot b_{i}$,
$R^{23}=\sum \id\ot a_{i}\ot b_{i}$ if $R$ has the form
$R=\sum a_{i}\ot b_{i}$.
We employ the following standard notation for the q-exponential:
\bn
\exp_{q}(x):=1+x+\mbox{\ls$\frac{x^{2}}{(2)_{q}!}$}+\ld+
\mbox{\ls$\frac{x^{n}}{(n)_{q}!}$}+\ld=\sum_{n\geq0}
\mbox{\ls$\frac{x^{n}}{(n)_{q}!}$}~,
\lb{UR3}
\ed
where
\bn
(n)_{q}:=\mbox{\ls$\frac{q^{n}-1}{q-1}$}~.
\lb{UR4}
\ed
A general formula for the universal $R$-matrix $R$ for quantized
contragredient Lie (super)algebras was presented in Refs. 
\cite{KT3}--\cite{KT5}. Here we specialize this result to our case 
$U_q(g)$, where $g=A_1^{(1)},~C(2)^{(2)}$.

The explicit expression of the universal $R$-matrix $R$ for our case 
$U_q(g)$ can be presented as follows:
\bn
R=R_{+}R_{0}R_{-}K~.
\lb{UR5}
\ed
Here the factors $K$ and $R_{\pm}$ have the following form
\bneqn
K&\!\!=\!\!&q^{\frac{1}{(\a,\a)}h_{\a}\ot h_{\a}+
h_{\d}\ot h_{\rm d}+h_{\rm d}\ot h_{\d}}~,
\lb{UR6}
\\[9pt]
R_{+}&\!\!=\!\!&\prod_{n\ge0}^{\ra}\cR_{n\d+\a}, \qquad
R_{-}=\prod_{n\ge1}^{\la}\cR_{n\d-\a}.
\lb{UR7}
\edeqn
The elements $\cR_{\g}$ are given by the formula
\bn
\cR_{\g}=\exp_{\bar{q}_{\g}}\Bigl(A(\g)
(q-q^{-1})(e_{\g}\ot e_{-\g})\Bigr)~,
\lb{UR8}
\ed
where
\bneqn
A(\g)&\!\!=\!\!&\left\{\begin{array}{lll}
&\!\!\!\!\!(-1)^{n\th(\a)}\quad &{\rm if}\;\g=n\d+\a~,\\
&\!\!\!\!\!(-1)^{(n-1)\th(\a)} \quad &{\rm if}\;\g=n\d-\a~.
\end{array}\right.
\lb{UR9}
\edeqn
Finally, the factor $R_{0}$ is defined as follows
\bn
R_{0}=\exp\Bigl((q-q^{-1})\sum_{n>0}^{}d(n)
e_{n\d}\ot e_{-n\d}\Bigr)~,
\lb{UR10}
\ed
where $d(n)$ is the inverse to $a(n)$, i.e.
\bn
d(n)=\mbox{\ls$\frac{n(q-q^{-1})}
{q^{n(\a,\a)}-q^{-n(\a,\a)}}$}~.
\lb{UR11}
\ed

\setcounter{equation}{0}
\section{The 'new realization'}
Let us denote by $d$ the Cartan element $h_{\rm d}$ and by $c$
the Cartan element $h_\d$, emphasizing that $d$ defines homogeneous
gradation of the algebra and $k_{\d}\!=\!q^{h_{\d}}$ is the central
element. It will be convenient in the following to add its square
roots $q^{\pm\frac{c}{2}}\!=\!k_{\d}^{\pm\frac{1}{2}}$.
Let us introduce the new notations:
$e_{n}^{}\!:=\!e_{n\d+\a}^{}$ ($n\geq0$),
$e_{-n}^{}\!:=\!-(\!-1)^{(n-1)\th(\a)}k_{-n\d+\a}^{}e_{-n\d+\a}^{}$
($n>0$), and $f_{n}^{}\!:=\!-e_{n\d-\a}^{}k_{n\d-\a}^{}$ ($n>0$),
$f_{-n}^{}\!:=\!(\!-1)^{(n+1)\th(\a)}\!e_{-n\d-\a}^{}$ ($n\geq0$).
We also put $a_{n}^{}\!:=\!e_{n\d}^{}q^{\frac{nc}{2}}$
($n\geq1$), and $a_{-n}^{}\!:=\!(\!-1)^{n\th(\a)}\!e_{-n\d}^{}
q^{-\frac{nc}{2}}$ ($n\geq1$).
Collect the elements $e_n$, $f_n$ ($n\in\ZZ\,$) and $a_{\pm n}$
($n\geq1$) into the generating functions ("fields")
\bn
\begin{array}{rcccl}
e(z)&\!\!=\!\!&\sum\limits_{n\in\,\ZZ}e_nz^{-n},\qquad
\psi_{+}(z)&\!\!=\!\!&k_\a^{-1}\exp\Bigl((\!-1)^{\th(\a)}(q-q^{-1})
\sum\limits_{n=1}^{\infty}a_{n}z^{-n}\Bigr),
\\[11pt]
f(z)&\!\!=\!\!&\sum\limits_{n\in\,\ZZ}f_nz^{-n},\qquad
\psi_{-}(z)&\!\!=\!\!&k_{\a}\exp\Bigl((\!-1)^{\th(\a)}(q^{-1}-q)
\sum\limits_{n=1}^{\infty}a_{-n}z^{n}\Bigr),
\lb{NR1}
\end{array}
\ed
such that
\bn
{\rm deg}\, e(z)\!=\!{\rm deg}\, f(z)\!=\!\th(\a),
\qquad {\rm deg}\,\psi_\pm(z)\!=\!0.
\ed
These fields satisfy the following conjugation conditions
with respect to graded conjugation
"$^{\ddagger}$"}:
\bn
\begin{array}{rcccl}
(e(z))^{\ddagger}&\!\!=\!\!&f(z^{-1})~,\qquad\quad\;
(f(z))^{\ddagger}&\!\!=\!\!&(\!-1)^{\th(\a)}e(z^{-1})~,
\\[7pt]
(\psi_{+}(z))^{\ddagger}&\!\!=\!\!&\psi_{-}(z^{-1})~,\qquad
(\psi_{-}(z))^{\ddagger}&\!\!=\!\!&\psi_{+}(z^{-1})~,
\end{array}
\lb{NB1'}
\ed
and have the following
symmetry with respect to the translation operator $T_{\d}$:
\bn
\begin{array}{rcccl}
T_{\d}(e(z))&\!\!=\!\!&(-1)^{\th(\a)}z^{}e((\!-1)^{\th(\a)}z)~,\quad\;
T_{\d}(f(z))&\!\!=\!\!&(-1)^{\th(\a)}z^{-1}f((\!-1)^{\th(\a)}z)~,
\\[7pt]
T_{\d}(\psi_{+}(z))&\!\!=\!\!&q^{-c}\psi_{+}((\!-1)^{\th(\a)}z)~,
\qquad\;\;
T_{\d}(\psi_{-}(z))&\!\!=\!\!&q^{c}\psi_{-}((\!-1)^{\th(\a)}z)~.
\end{array}
\lb{NB1"}
\ed
\bnpr
In terms of the fields (\ref{NR1}) the relations of Section 4
can be rewritten in the following compact form
\bn
\begin{array}{rcl}
[q^c,{\rm everything}]\!\!&=\!\!&0~,
\\[7pt]
u^d\vph(v)u^{-d}\!\!&=\!\!&\vph(uv)\qquad
\lb{NR2}
\end{array}
\ed
where $\vph(v)=e(v),f(v),\psi_{\pm}(v)$, and also
\bneqn
\psi_\pm(u)\psi_\pm(v)\!\!&=\!\!&\psi_\pm(v)\psi_\pm(u)~,
\lb{NR3}
\\[4pt]
(u-\bar{q}_{\a}v)e(u)e(v)\!\!&=\!\!&(\bar{q}_{\a}u-v)e(v)e(u)~,
\lb{NR4}
\\[4pt]
(u-q_{\a}v)f(u)f(v)\!\!&=\!\!&(q_{\a}u-v)f(v)f(u)~,
\lb{NR5}
\\[4pt]
\psi_\pm(u)e(v)\bigl(\psi_\pm(u)\bigr)^{-1}\!\!&=\!\!&(\!-1)^{\th(\a)}
\mbox{\ls$\frac{\bar{q}_{\a}q^{\mp\frac{c}{2}}u-v}
{q^{\mp\frac{c}{2}}u-\bar{q}_{\a}v}$}\,e(v)~,
\lb{NR6}
\\[4pt]
\psi_\pm(u)f(v)\bigl(\psi_\pm(u)\bigr)^{-1}\!\!&=\!\!&(\!-1)^{\th(\a)}
\mbox{\ls$\frac{q_{\a}q^{\pm\frac{c}{2}}u-v}
{q^{\pm\frac{c}{2}}u-q_{\a}v}$}\,f(v)~,
\lb{NR7}
\\[4pt]
\bigl(\psi_+(u)\bigr)^{-1}\psi_{-}(v)\psi_+(u)
\bigl(\psi_-(v)\bigr)^{-1}\!\!&=\!\!&
\mbox{\ls$\frac{(q^cu-q_{\a}v)(q^{-c}u-\bar{q}_{\a}v)}
{(q^cv-\bar{q}_{\a}u)(q^{-c}v-q_{\a}u)}$}~,
\lb{NR8}
\edeqn
\bn
[e(u),f(v)]\!=\!\mbox{\ls$\frac{1}{q-q^{-1}}$}
\Bigl(\d(\mbox{\ls$\frac{u}{v}$}q^{-c})\psi_-(vq^{c/2})-
\d(\mbox{\ls$\frac{u}{v}$}q^{c})\psi_+(uq^{c/2})\Bigr)~.\qquad\quad
\lb{NR9}
\ed
\lb{PNR1}
\edpr
Here in (\ref{NR8}) $\d(z)=\sum_{n\in\,\ZZ}z^n$, and the brackets
$[\cdot, \cdot]$ in the relation (\ref{NR9}) mean the supercommutator:
\bn
[e(u),f(v)]=e(u)f(v)-(-1)^{\th(\a)}f(v)e(u)~.
\lb{NR10}
\ed
Given description is called the 'new realization', or the current
realization of the quantum affine superalgebras $U_q(A_1^{(1)})$
and $U_q(C(2)^{(2)})$. It should be noted that the relations
(\ref{NR1}) and (\ref{NR3})--(\ref{NR9}) differ from the corresponding
relations of Refs. \cite{DK}, \cite{D}, \cite{YZ} by replacement
$q$ to $q^{-1}$.
The current realization possesses its own graded comultiplication
structure, different from \r{DR15}:
\bn
\begin{array}{rcl}
\DD(c)\!\!&=\!\!&c\ot1+1\ot c~,
\\[6pt]
\DD(d)\!\!&=\!\!&d\ot1+1\ot d~,
\\[6pt]
\phantom{pos}\qquad\DD(\psi_\pm(z))\!\!&=\!\!&
\psi_\pm(zq^{\pm\frac{c_2}{2}})\ot\psi_\pm(zq^{\mp\frac{c_1}{2}})~,
\\[6pt]
\DD(e(z))\!\!&=\!\!&e(z)\ot 1+
\psi_-(zq^{\frac{c_1}{2}})\ot e(zq^{c_1})~,
\\[6pt]
\DD(f(z))\!\!&=\!\!&f(zq^{c_2})\ot\psi_+(zq^{\frac{c_2}{2}})+
1\ot f(z)~,
\lb{NR11}
\end{array}
\ed

\bn
\begin{array}{rcl}
\SD(c)\!\!&=\!\!&-c~,\quad \SD(d)=-d~,
\\[6pt]
\SD(\psi_\pm(z))\!\!&=\!\!&\bigl(\psi_\pm(z)\bigr)^{-1}~,\;\;
\\[6pt]
\SD(e(z))\!\!&=\!\!&-\bigl(\psi_{-}(zq^{-\frac{c}{2}})\bigr)^{-1}
e(zq^{-c})~,
\\[6pt]
\SD(f(z)\!\!&=\!\!&-f(zq^{-c})
\bigl(\psi_{+}(zq^{-\frac{c}{2}})\bigr)^{-1}~,
\lb{NR12}
\end{array}
\ed

\bn
\vep(c)=\vep(d)=\vep(e(z))=\vep(f(z))=0~,\qquad
\vep(\psi^\pm(z))=1~.\quad
\lb{NR13}
\ed
Here $\DD$, $\SD$, and $\vep$ are the comultiplication, antipode
and counite correspondingly.
The two comultiplications $\D_q$ and $\DD$ are related by the twist
\cite{KT5}:
\bn
\DD(x)=F^{-1}\D(x)F~,
\lb{NR14}
\ed
where $F=R_+^{21}$, with $R_+$ given by (\ref{UR7})--(\ref{UR9}), 
such that the universal $R$-matrix for the comultiplication $\DD$ 
equals to
\bn
\RD=R_0R_-KR_+^{21}
\lb{NR15}
\ed
with the factors from (\ref{UR5}). In the generators $e_n$, $f_n$ and
$a_n$ it can be rewritten as follows
\bn
\RD=\cK\kr~,\qquad
\lb{NR16}
\ed
where
\bneqn
\cK\!\!&=\!\!&q^{\frac{h_{\a}\ot h_{\a}}{(\a,\a)}}
q^{\frac{1}{2}(c\ot d+d\ot c)}\exp\Bigl((q-q^{-1})
\sum_{n=1}^{\infty}\!\ol{d}(n)\,a_{n}\ot a_{-n}\Bigr)
q^{\frac{1}{2}(c\ot d+d\ot c)}~,
\lb{NR17}
\\[5pt]
\kr\!\!&=\!\!&\prod_{n\in\ZZ}^{\ra}
\exp_{\bar{q}_{\a}}\!\Big((q^{-1}-q)f_{-n}\ot e_{n}\Big)~,
\lb{NR18}
\lb{PNR2}
\edeqn
and
\bn
\ol{d}(n)=\mbox{\ls$\frac{n(q-q^{-1})}{q_\a^{n}-q_\a^{-n}}$}~.
\ed
It is possible to give another presentation of the element $\kr$ in
the completed algebras $\bar{U}(g)$, where $g$ is either $A_1^{(1)}$
or $C(2)^{(2)}$ \cite{DK}, \cite{DKP}. The completion is done with
respect to open neighborhoods of zero $\bar{U}_r=\sum _{s>r}U_s$,
where $U_s$ consists of all the elements from $U(g)$ of degree $s$.
The completed algebra acts on  (infinite-dimensional) representations
of highest weight and admits the series over monomials
$x_{i_1}x_{i_2}\cd x_{i_n}$, $i_1\leq i_2...\leq i_n,$ with $x=e,f,a$
and fixed $\sum i_k$. The matrix coefficients of the products of the
currents $e(z_1)e(z_2)\cd e(z_n)$ and $f(z_1)f(z_2)\cd f(z_n)$,
defined originally as formal series, converge to meromorphic in
${\CC}^n$ functions with the poles at $z_i=0$ and
$z_i=q_{\a}^{\mp 1}z_j$, $i\leq j$.

Let $t(z)=(q-q^{-1})f(z)\ot e(z)$. As before, we understand the
product $t(z_1)\cd t(z_n)$ as operator-valued meromorphic function
in $\bigl(\CC^*\bigr)^n$ with simple poles at $z_i=q_{\a}^{\mp1}z_j$,
$i\not=j$.
Define
\bn
{\kr}'=1+\sum_{n>0}\mbox{\ls$\frac{1}{n!(2\pi i)^n}$}
\stackreb{D_n}{\oint\cd\oint}\mbox{\ls$\frac{dz_1}{z_1}$}\cd
\mbox{\ls$\frac{dz_n}{z_{n}}$}\,t(z_1)\cd t(z_n)~,
\lb{NR19}
\ed
and integration region $D_n$ is defined as
$D_n=\Bigl\{\Bigl|z_i|\!=\!1,\;i\!=\!1,...,n\Bigr\}\;$ for $|q|<1$
and, more generally, by
\bn
D_n=\Bigl\{\Bigl|z_i\!\!\!\!\prod_{j=1,...,n,\atop j\not=i}
\!\!\!(z_i-q_{\a}z_j)\Bigr|\!=\!1, \; i=1,...,n\Bigr\}
\lb{NR20}
\ed
for any $q$, such that $q_{\a}^N\not=1$, $N\in\ZZ\setminus\!\{0\}$.
\bnpr
The action of the tensor ${\cR}'={\cK}\kr'$ in tensor product of
highest weight modules is well defined and coincides with the action
of the universal $R$-matrix \r{NR16}
\edpr
The integrals in \r{NR17} can be computed explicitly.
Let us put by induction
\bn
t^{(n)}(z)=-\!\!\!\!\!\!
{\mathop{\rm Res}\limits_{\quad z_1=z\bar{q}_{\a}^{2n-2}}}
\!t(z_1)t^{(n-1)}(z)\mbox{\ls$\frac{dz_1}{z_1}$},
\lb{NR21}
\ed
where $t^{(1)}(z)=t(z)$. In the components the fields $t^{(n)}(z)$
look as follows:
\bn
t^{(n)}(z)=C_n e^{(n)}(z)\ot f^{(n)}(z)~,
\lb{NR22}
\ed
where
\bn
\begin{array}{rcl}
C_n\!\!&=\!\!&(-1)^{(n-1)\th(\a)}(q\!-\!q^{-1})^n
\tl{q}_{\a}^{-\frac{n(n-1)}{2}}(\tl{q}_{\a}\!-\!1)^{n-1}
(n\!-\!1)_{\tl{q}_{\a}}\!!\,(n)_{\tl{q}_{\a}}\!!~,
\\[9pt]
e^{(n)}(z)\!\!&=\!\!&e(z)e(\bar{q}_{\a}^{}z)e(\bar{q}_{\a}^{\;2}z)
\cd e(\bar{q}_{\a}^{\;n-2}z)e(\bar{q}_{\a}^{\;n-1}z)~,
\\[9pt]
f^{(n)}(z)\!\!&=\!\!&f(\bar{q}_{\a}^{\;n-1}z)
f(\bar{q}_{\a}^{\;n-2}z)\cd f(\bar{q}_{\a}^{\;2}z)
f(\bar{q}_{\a}z)f(z)~,
\end{array}
\lb{NR23}
\ed
such that
\bn
\begin{array}{rcl}
{e}^{(n)}(z)\!\!&=\!\!&\sum\limits_{m\in\,\ZZ}
\bigl(z\bar{q}_{\a}^{n}\bigr)^m
\!\!\sum\limits_{{\lm_1\geq \cd\geq\lm_n,\atop\lm_1+...+\lm_n=m}}
\mbox{\ls$\frac{q_{\a}^{\lm_1+2\lm_2+...n\lm_n}}
{\prod\limits_{j\in\,\ZZ}(\lm'_j-\lm'_{j+1})_{\bar{q}_{\a}}!}$}\,
e_{\lm_n}^{}e_{\lm_{n-1}}^{}\cd e_{\lm_1}^{}~,
\\[20pt]
{f}^{(n)}(z)\!\!&=\!\!&\sum\limits_{m\in\,\ZZ}\left(zq_{\a}^{}\right)^m
\!\!\sum\limits_{\lm_1\geq \cd\geq\lm_n,\atop\lm_1+...+\lm_n=m}
\mbox{\ls$\frac{q_{\a}^{\lm_1+2\lm_2+...n\lm_n}}
{\prod\limits_{j\in\,\ZZ}(\lm'_j-\lm'_{j+1})_{q_{\a}^{}}!}$}\,
f_{\lm_n}^{}f_{\lm_{n-1}}^{}\cd f_{\lm_1}^{}~.
\lb{NR24}
\end{array}
\ed
Here $\lm'_{j}\!=\!\#k$, such that $\lm_k\geq j$, and
$\tl{q}_{\a}\!:=\!q^{(\a,\a)}$. The product in denominator
is finite, since there are only finitely many distinct $\lm'_j$ for
a given choice of $\lm_k$. Then, repeating the calculations in \cite{DKP},
we get vertex type presentation of the element $\kr'$:
\bn
\kr'=\exp\Bigl(\sum_{n>0}\mbox{\ls$\frac{1}{n}$}I_n\Bigr)~,
\lb{NR25}
\ed
where the sequence of operators
\bn
I_n=\oint\mbox{\ls$\frac{t^{(n)}(z)dz}{2\pi iz}$}
\lb{NR26}
\ed
commute between themselves:
$$
[I_n,I_m]=0,\qquad n,m>0.
$$
The vertex operator presentation (\ref{NR25}) is convenient
for applications to integrable representations: it is expressed
through integrals over the fields, which number is precisely $k$
for level $k$ integrable representations.

\section{Final Remarks}
The aim of this paper is to describe in unified way with detail
the $q$-deformed untwisted affine algebra
$U_{q}(\widehat{sl}(2))\!=\!U_{q}(A_{1}^{(1)})$ and twisted
superalgebra $U_{q}(osp(2|2))^{(2)}\!=\!U_{q}(C(2)^{(2)})$.
In order to describe the complete list of quantum affine
(super)algebras of rank 2 one should consider 
some more three quantum affine (super)algebras:
$U_{q}(sl(1|3))^{(4)}\!=\!U_{q}(A{(0,2)})^{(4)}$,
$U_{q}(sl(3))^{(2)}\!=\!U_{q}(A_{2}^{(2)})$ and
$U_{q}(\widehat{osp}(1|2))\!=\!U_{q}(B{(0,1)}^{(1)}$.
The Dynkin diagram of the superalgebra $A(0,2)^{(4)}$ has as
geometric structure as the (super)algebras $A_{1}^{(1)}$ and
$C(2)^{(2)}$ but in this case the root $\a$ is even and $\d-\a$
is odd one, and the sector of imaginary roots has odd roots.
Therefore in the case of the quantum superalgebra $A(0,2)^{(4)}$
the relations of the type (\ref{CW28})--({\ref{CW31}) are
more complicated and they demand special consideration.
The second family of two quantum affine (super)algebras
$U_{q}(A_{2}^{(2)})$ and $U_{q}(B{(0,1)}^{(1)})$ are
described by the same Dynkin diagram with different colors of roots.
Preliminary results in this direction are given in \cite{LT}, where
in particular the Cartan-Weyl basis of basic affine superalgebra
$U_{q}(\widehat{osp}(1|2))$ is considered. The unified description
of three mentioned above quantum affine (super)algebras, analogous to
the one given in the present paper, is in preparation.

\section*{Acknowledgments}
This work was supported (S.M. Khoroshkin, V.N. Tolstoy) by the
Russian Foundation for Fundamental Research, grant No.98-01-00303,
by the program of French-Russian scientific cooperation
(CNRS grant PICS-608 and grant RFBR-98-01-22033), as well as by
KBN grant 2P03B13012 (J. Lukierski) and INTAS-99-1705 (S. Khoroshkin).

\end{document}